\pdfoutput=1
\documentclass[11pt,a4paper]{article}
\usepackage[utf8]{inputenc}
\usepackage{amsfonts, amsmath, amssymb, color, a4wide}
\usepackage[utf8]{inputenc} % allow utf-8 input
\usepackage[T1]{fontenc}    % use 8-bit T1 fonts
\usepackage{url}            % simple URL typesetting
\usepackage{booktabs}       % professional-quality tables
\usepackage{amsfonts}       % blackboard math symbols
\usepackage{nicefrac}       % compact symbols for 1/2, etc.
\usepackage{microtype}      % microtypography

\usepackage{lmodern}

\usepackage{amsfonts}
\usepackage{amsmath}
\usepackage{amssymb}
\usepackage{bm}
\usepackage{algorithm}
\usepackage{algorithmic}
\usepackage{mathtools}
\usepackage{xspace}
\usepackage{xcolor}
\usepackage{amsthm}
\usepackage{cases}
\usepackage{natbib}
\usepackage{bbm}
\usepackage{dsfont}
\RequirePackage[shortlabels]{enumitem}

\newtheorem{theorem}{Theorem}[section]

\newtheorem{lemma}[theorem]{Lemma}
\newtheorem{corollary}[theorem]{Corollary}

\newtheorem*{remark}{Remarks}

\newtheorem{assumption}{Assumption}

\input{macros_gilles_2.tex}

\DeclareMathOperator{\Var}{Var}

\newenvironment{sproof}{%
	\proof}{\endproof}

\makeatletter
\newenvironment{protocol}[1][htb]{%
	\renewcommand{\ALG@name}{Protocol}% Update algorithm name
	\begin{algorithm}[#1]%
	}{\end{algorithm}}
\makeatother

\newcommand{\vertiii}[1]{{\left\vert\kern-0.25ex\left\vert\kern-0.25ex\left\vert #1 
		\right\vert\kern-0.25ex\right\vert\kern-0.25ex\right\vert}}

\title{Constant regret for sequence prediction with limited advice}
\author{El Mehdi Saad$^{1}$ , Gilles Blanchard$^{1,2}$ \\
	$^1$Laboratoire de Mathématiques d'Orsay, CNRS, Université Paris-Saclay; $^2$Inria\\}
\date{}

\begin{document}
	\maketitle

\begin{abstract}
	We investigate the problem of cumulative regret minimization for individual sequence prediction 
	with respect to the best expert in a finite family of size $K$ under limited access to information. We assume that in each round, the learner can predict using
	a convex combination of at most $p$ experts for prediction, then they can observe {\em a posteriori} the losses of at most $m$ experts. We assume that the loss function is range-bounded and exp-concave. In the standard multi-armed bandits setting, when the learner is allowed to play only one expert per round and observe only its feedback, known optimal regret bounds are of the order $\cO(\sqrt{KT})$. We show that allowing the learner to play one additional expert per round and observe one additional feedback improves substantially the guarantees on regret. We provide a strategy combining only $p=2$ experts per round for prediction and observing $m\ge 2$ experts' losses.
	Its randomized regret
	(wrt. internal
	randomization of the learners' strategy) is of order $\mathcal{O}\left((K/m)\log(K\delta^{-1})\right)$ with probability $1-\delta$,
	i.e., is independent of the horizon $T$ (``constant'' or ``fast rate'' regret) if ($p\ge 2$ and $m\ge 3$). We prove that this rate is optimal up to a logarithmic factor in $K$. In the case $p=m=2$, we provide an upper bound of order $\mathcal{O}(K^2\log(K\delta^{-1}))$, with probability $1-\delta$. Our strategies do not require any prior knowledge of the horizon $T$ nor of the confidence parameter $\delta$. Finally, we show that if the learner is constrained to observe only one expert feedback per round, the worst-case regret is the ``slow rate'' $\Omega(\sqrt{KT})$, suggesting that synchronous observation of at least two experts per round is necessary to have a constant regret. 
\end{abstract}


\section{Introduction}

We study the problem of online individual sequence prediction with expert advice, based on the setting presented by \citet[Chap. 2]{cesa2006prediction}, under limited access to information. In this game, the learner's aim is to predict an unknown sequence $(y_1,y_2,\dots)$ of an outcome space $\mathcal{Y}$. The mismatch between the learner's predictions $(z_1,z_2,\dots )$, taking values in a closed convex subset $\mathcal{X}$ of a real vector space, and the target sequence is measured via a loss function $\ell(z,y)$.
The learner's predictions may only depend on past observations. Following standard terminology
used in prediction games, we will use the word ``play'' to mean the prediction output by the learner.

In each round $t \in \intr{T}$ (for a non-negative integer $n$, we denote $\intr{n} = \{1, \dots, n\}$), the learner has access to $K$ experts predictions $(F_{1,t}, \dots, F_{K,t})$. The performance of the learner is compared to that of the best single expert. More precisely, the objective is to have a cumulated regret as small as possible, where the regret is defined by
\begin{equation*}
\mathcal{R}_T =  \sum_{t=1}^{T} \ell\left(z_t, y_t \right) - \min_{i\in \intr{K}} \sum_{t=1}^{T}\ell\left(F_{i,t}, y_t\right).
\end{equation*}

\noindent Experts aggregation is a standard problem in machine learning, where the learner observes the predictions of all experts in each round and plays a convex combination of those. However,
in many practical situations, querying the advice of every expert is unrealistic. Natural constraints arise, such as the financial cost of consultancy, time limitations in online systems, or computational
budget constraints if each expert is actually the output of a complex prediction model. One might hope to make predictions in these scenarios while minimizing the underlying cost. Furthermore, we will distinguish between the constraint on the number of experts' advices used for prediction, and the number of feedbacks (losses of individual experts) observed {\em a posteriori}. This difference naturally arises in online settings where the advices are costly prior to the prediction task but just observing reported
experts' losses after prediction can be cheaper.
If the learner picks one single expert per round, plays the prediction of that expert and observes the
resulting loss, the game is the standard multi-armed bandits problem. In this paper, we investigate intermediate settings, where the player has a constraint $p \le K$ on the number of experts used for prediction (via
convex combination) in each round and several feedbacks $m\le K$ of actively chosen experts to see their losses. In the standard multi-armed bandit problem, the played arm is necessarily the observed arm, this restriction is known as the \emph{coupling between exploitation and exploration}. In our protocol, we consider a generalization of that restriction through the \emph{Inclusion Condition (IC)}: when $m\ge p$, if $\text{IC}= \text{True}$, we require that the set of played experts for prediction at round $t$, denoted $S_t$ , is included in the set of observed experts, denoted $C_t$. More precisely, if $\text{IC}= \text{True}$, in each round $t$, the player first chooses $p$ experts out of $K$ and plays a convex combination of their prediction, then she observes the feedback (loss) of the individual selected experts, then picks $m-p$ additional experts to observe their losses. When $\text{IC} = \text{False}$, the choice of played and observed experts is decoupled; this means that the loss incurred by the $p$ experts used
for prediction is not necessarily observed.

\begin{protocol}
	%\centering
	\caption{ The Game Protocol $(p,m,\text{IC})$.}
	\begin{algorithmic}\label{algo:gp}
		\STATE \textbf{Parameters}:
		\STATE \quad $p$, the number of experts allowed for prediction.
		\STATE \quad $m$, the number of experts allowed for observation as feedback.
		\STATE \quad $\bm{\text{IC}} \in \{\text{False},\text{True}\}$, inclusion condition (if $\bm{\text{IC}} = \text{True}$, we must have $p\le m$).
		\FOR{\textbf{each} round $t=1,2,\dots,T$}
		\STATE Choose a subset $S_t \subseteq \intr{K}$ such that $|S_t| = p$,
		and convex combination weights $(\alpha_i)_{i \in S_t}$.
		\STATE Play the convex combination $\sum_{i \in S_t} \alpha_{i,t} F_{i,t}$ and incur its loss.
		\IF{$\bm{\text{IC}} = \text{True}$,}
		\STATE Choose a subset $C_t \subseteq \intr{K}$ such that: $\lvert C_t \rvert = m$ and $S_t\subseteq C_t$.
		\ELSIF{$\bm{\text{IC}} = \text{False}$,} 
		\STATE Choose a subset $C_t \subseteq \intr{K}$ such that: $\lvert C_t \rvert = m$.
		\ENDIF
		\STATE The environment reveals the losses $\left(\ell\left(F_{i,t},y_t\right)\right)_{i\in  C_t}$.
		\ENDFOR
	\end{algorithmic}
\end{protocol}

% \begin{remark}
% 	The case $p=1$ and $m=1$ reduces to the classical multi-armed bandit setting (\citealp{bubeck2012regret}), where the played arm is necessarily the observed arm in each round. The learner is faced with the exploration-exploitation dilemma referring to  balancing the need to query various arms and sticking with the best-performing arm so far. The condition $\textbf{IC}=\text{True}$ can be seen as a generalization of the coupling between exploration and exploitation through the restriction: $S_t\subseteq C_t$. When $\textbf{IC}=\text{False}$, we consider that the played and observed experts are decoupled (they can be freely chosen by the learner). The decoupled setting was motivated and analysed in \cite{avner2012decoupling} for $p=m=1$. 
% \end{remark}

A closely related question was considered by~\cite{seldin2014prediction}, obtaining $\cO(\sqrt{T})$
regret bounds for a general loss function (see extended discussion in the next section.)
Our emphasis here is on obtaining \textit{constant bounds} guarantees on regret (i.e. independent of the time horizon $T$). Such ``fast" rates, linked to assumptions related to strong convexity of the loss function $\ell$,
%, are considered important and
have been the subject of many works in learning (batch and online, in the stochastic setting) and optimization, but are comparatively under-explored in fixed sequence prediction.

In the literature on the prediction of fixed individual sequences, no assumptions are made about the distribution of the sequences. 
The attainability of fast rates (or constant regrets) is also possible under certain assumptions
on the loss function $\ell$:
the full information setting was studied, mainly by \cite{vovk1990aggregating}, \cite{vovk1998game}, \cite{vovk2001competitive}, where it was shown that fast rates are attainable under the $\textit{mixability}$ assumption on the loss function. The reader can find an extensive discussion of different assumptions considered in the literature for this problem in \cite{van2015fast}.   
%Below we present the only
In the present paper, we make the following assumption on the loss function:
\pagebreak
%made in this paper , namely range-boundedness and exp-concavity.
\begin{assumption}\label{assump}
	There exist $B, \eta >0$, such that
	\begin{itemize}
		\item \textbf{Exp-concavity:} For all $y \in \mathcal{Y}$, $\ell(.,y)$ is $\eta$-exp-concave over domain $\cX$.
		\item \textbf{Range-boundedness:} For all $y \in \mathcal{Y}$: $\quad \sup_{x,x' \in \cX} \abs{\ell(x,y) - \ell(x',y)} \le B$.
	\end{itemize}
\end{assumption}

\begin{remark}
	This assumption is satisfied in some usual settings of learning theory such as the least squares loss with bounded outputs: $\mathcal{X}=\mathcal{Y} = [x_{\min}, x_{\max}]$ and $\ell(x, x') = (x-x')^2$. Then $\ell$ satisfies Assumption~\ref{assump}, with $B =  (x_{\max} - x_{\min})^2$ and $\eta = 1/(2B)$.
	% The second assumption is implied by: $\ell(x,y) \in [0,B]$. Note that
\end{remark}

\begin{remark}
	The regret as well as all the algorithms to follow remain unchanged if we replace $\ell$ by $\tilde{\ell}:\cX \to [0,B]$ defined by $\tilde{\ell}(x,y) := \ell(x,y) - \min_{x \in \cX}\ell(x,y)$, so we can assume without loss of generality $\ell \in [0,B]$ instead of range-boundedness;
	the results obtained still hold in the latter more general case.
\end{remark}

Assumption~\ref{assump} was considered in several previous works tracking fast rates both in batch and online learning (\citealp{koren2015fast}, \citealp{mehta2017fast}, \citealp{gonen2016tightening}, \citealp{mahdavi2015lower}, \citealp{van2015fast}). We introduce a new characterization for the class of functions satisfying Assumption~\ref{assump}. Let $c>0$, define $\mathcal{E}(c)$ as the class of functions $f: \mathcal{X} \to \mathbb{R}$, such that
\begin{equation}\label{eq:prop}
\forall x, x' \in \mathcal{X}:\quad f\left(\frac{x+x'}{2}\right) \le \frac{1}{2} f(x) + \frac{1}{2} f(x') - \frac{1}{2c} \left(f(x) - f(x')\right)^2.
\end{equation}
%The class of functions $\cE(.)$ was i
We introduce this class to highlight the sufficient and minimal property of $\ell$ required for the
proofs in this paper to work, namely we will only make use of~\eqref{eq:prop} in the proofs
of the results to come.
%namely: property \eqref{eq:prop}. 
%\end{remark}

Lemma~\ref{lem:assumption1} below relates the class of functions $\cE(.)$ to the set of functions satisfying Assumption~\ref{assump} as well a sufficient condition (Lipschitz and Strongly Convex or LIST condition). \pagebreak
\begin{lemma}\label{lem:assumption1}
	Let $y \in \mathcal{Y}$ be fixed.
	\begin{itemize}
		\item If $\ell(.,y)$ is $B$-range-bounded and $\eta$-exp-concave, then: $ \ell(.,y) \in \mathcal{E}\left( \frac{\eta B^2}{4 \log\left(1 + \frac{\eta^2B^2}{2}\right)}\right)$.
		\item If $\ell(.,y) \in \mathcal{E}(c)$ and is continuous, then: $\ell(.,y)$ is $c$-range-bounded and $(4/c)$-exp-concave.
		\item If $\ell(.,y)$ is $L$-Lipschitz and $\rho$-strongly convex, then $\ell(.,y) \in \cE(4L^2/\rho)$.
	\end{itemize}
\end{lemma}

Figure~\ref{fig:table} summarizes bounds on regret for bounded and exp-concave loss functions. We only consider fixed individual sequences, which corresponds to fully oblivious adversaries (see \citealp{audibert2010regret} for a definition of different types of adversaries).
\begin{figure}[!ht]\label{fig:table}
	\begin{center}
		\begin{tabular}{ c|c|c|c|c| } 
			\cline{2-5}
			{} & \multicolumn{2}{c|}{} & \multicolumn{2}{c|}{}  \\ 
			{} & \multicolumn{2}{c|}{$p=1$} & \multicolumn{2}{c|}{$p\ge2$}  \\ 
			{} & \multicolumn{2}{c|}{} & \multicolumn{2}{c|}{}  \\  \hline
			\multicolumn{1}{|c|}{}& {Lower bound} & {Upper bound} & {Lower bound} & {Upper bound $(p=2)$} \\ 
			\cline{2-5}
			\multicolumn{1}{|c|}{} & {} & {} & {} & {}   \\ 
			\multicolumn{1}{|c|}{$m=1$} & $\sqrt{KT}$ & $\sqrt{KT}$ & $\sqrt{KT}$ & $\sqrt{KT}$  \\ 
			\multicolumn{1}{|c|}{} & {[1]} & {[2]} & {[Thm~\ref{thm:4}]} & {[2]} \\
			\multicolumn{1}{|c|}{} & {} & {} & {} & {} \\ \hline
			\multicolumn{1}{|c|}{} & {} & {} & {} & {} \\
			\multicolumn{1}{|c|}{} & {} & {} & {} & $\textbf{IC}=\text{True}:K^2\log(K)$ \\
			\multicolumn{1}{|c|}{$m=2$} & $\sqrt{KT}$ & $\sqrt{KT}$ & $K$ & $\textbf{IC}=\text{False}:K\log(K)$  \\ 
			\multicolumn{1}{|c|}{} & {[3]} & {[2]} & {[Thm~\ref{thm:2}]} & {[Thm~\ref{thm:1bis} and~\ref{thm:1}]} \\
			\multicolumn{1}{|c|}{} & {} & {} & {} & {} \\ \hline
			\multicolumn{1}{|c|}{} & {} & {} & {} & {} \\
			\multicolumn{1}{|c|}{$m \ge 3$} & $\sqrt{\frac{K}{m}T}$ & $\sqrt{\frac{K}{m}T\log(K)}$ & $\frac{K}{m}$ & $\frac{K}{m}\log(K)$  \\
			\multicolumn{1}{|c|}{} & {[3]} & {[3]} & {[Thm~\ref{thm:2}]} & {[Thm~\ref{thm:1}]} \\
			\multicolumn{1}{|c|}{} & {} & {} & {} & {} \\
			\hline
		\end{tabular}
	\end{center}
	\caption{Existing bounds from the literature ([1] = \citealp{auer2002nonstochastic}, [2]=\citealp{audibert2010regret}, [3]=\citealp{seldin2014prediction}) and new bounds presented in this paper. All bounds hold up to numerical constant factors. Under Assumption~\ref{assump}, all new upper bounds hold with high probability if we replace the factor $\log(K)$ with $\log(K\delta^{-1})$, $\delta$ being the confidence parameter. Lower bounds are in expectation. When bounds are the same, we omit the distinction between the settings $\textbf{IC}=\text{True}$ and $\textbf{IC}=\text{False}$ (coupling between exploration and exploitation, see Protocol~\ref{algo:gp}). }
\end{figure}

The remainder of this paper is organized as follows. Section~\ref{sec:related_w} presents some results from the literature relevant to the studied problem. Section~\ref{sec:Main1} introduces algorithms satisfying constant regrets in expectation in the case $p=2$ and $m\ge 3$; that section aims to present a preliminary view of the intuitions for attaining our objective. Next, we present in Section~\ref{sec:Main2} our main results consisting of algorithms satisfying constant regrets with a high probability for $p,m \ge 2$. %The latter section aims to present a preliminary view of the intuitions for attaining our objective.
Finally, in Section~\ref{sec:lower}, we present lower bounds for all the possible settings.

\section{Discussion of related work}\label{sec:related_w}

\paragraph{Games with limited feedback  and $\mathcal{O}\paren[1]{\sqrt{T}}$ regret:} In the standard setting of multi-armed bandit problem, the learner has to repeatedly obtain rewards (or incur losses) by choosing from a fixed set of $k$ actions and gets to see only the reward of the chosen action. Algorithms such as EXP3-IX \citep{neu2015explore} or EXP3.P \citep{auer2002nonstochastic} achieve the optimal regret of order $\mathcal{O}\paren[1]{\sqrt{KT}}$ up to a logarithmic factor, with high probability. A more general setting closer to ours was introduced by \citet{seldin2014prediction}. Given a budget $m\in \intr{K}$, in each round $t$, the learner plays the prediction of one expert $I_t$, then gets to choose a subset of experts $C_t$ such that $I_t \in C_t$ in order to see their prediction. A careful adaptation of the EXP3 algorithm to this setting leads to an expected regret of order $\mathcal{O}\paren[1]{\sqrt{(K/m)T}}$, which is optimal up to logarithmic factor in $K$.  

There are two significant differences between our framework and the setting presented by \cite{seldin2014prediction}. First, we allow the player to combine up to $p$ experts out of $K$ in each round for prediction. Second, we make an additional  exp-concavity-type assumption (Assumption~\ref{assump}) on the loss function. These two differences allow us to achieve constant regrets bounds (independent of $T$).

Playing multiple arms per round was considered in the literature of multiple-play multi-armed bandits. This problem was investigated under a budget constraint $C$ by \cite{zhou2018budget} and \cite{xia2016budgeted}. In each round, the player picks $m$ out of $K$ arms, incurs the sum of their losses. In addition to observing the losses of the played arms, the learner learns a vector of costs which has to be covered by a pre-defined budget $C$. Once the budget  is consumed, the game finishes. An extension of the EXP3 algorithm allows deriving a strategy in the adversarial setting with regret of order $\mathcal{O}\paren[1]{\sqrt{KC\log(K/m)}}$. The cost of each arm is supposed to be in an interval $[c_{\min}, 1]$, for a positive constant $c_{\min}$. Hence the total number of rounds in this game $T$ satisfies  $T = \Theta (C/m)$. Another online problem aims at minimizing the cumulative regret in an adversarial setting with a small effective range of losses. \cite{gerchinovitz2016refined} have shown the impossibility of regret scaling with the effective range of losses in the bandit setting, while \cite{thune2018adaptation} showed that it is possible to circumvent this impossibility result if the player is allowed one additional observation per round. However, it is impossible to achieve a regret dependence on $T$ better than the rate of order $\mathcal{O}\paren[1]{\sqrt{T}}$ in this setting.

Decoupling exploration and exploitation was considered by \cite{avner2012decoupling}. In each round, the player plays one arm, then chooses one arm out of $K$ to see its prediction (not necessarily the played arm as in the canonical multi-armed bandits problem). They devised algorithms for this setting and showed that the dependence on the number of arms $K$ can be improved. However, it is impossible to achieve a regret dependence on $T$ better than $\mathcal{O}\paren[1]{\sqrt{T}}$.

Prediction with limited expert advice was also investigated by \cite{helmbold1997some},\citet[Chap. 6]{cesa2006prediction} and \cite{cesa2005minimizing}. However, in these problems, known as label efficient prediction, the forecaster has full access to the experts advice but limited information about the past outcomes of the sequence to be predicted. More precisely, the outcome $y_t$ is not necessarily revealed to the learner. In such a framework, the optimal regret is of order $\mathcal{O}\paren[1]{\sqrt{T}}$.

\paragraph{Constant regrets in the full information setting:} 
The setting where the learner plays a combination of all the experts and is allowed to see all their predictions in each round is known in the literature as experts aggregation problem. It is a well-established framework \citep{cesa2006prediction} studied earlier by \cite{freund1997decision},  \cite{kivinen1999averaging}, \cite{vovk1998game}. This setting was investigated under the assumption that the loss $\ell$ function is $\eta$-exp-concave (i.e., the function $\exp(-\eta \ell)$ is concave). The Weighted Average Algorithm algorithm \citep{kivinen1999averaging} is known to achieve a constant regret of order $\mathcal{O}\left(\log(K)/\eta\right)$. While this result holds for any sequence of target variable and experts, it requires using a combination of all the experts in each round. In several situations, it is desirable to query and use the least number possible of experts advice for various reasons (such as cost or time restrictions). In this paper, we aim at achieving the same bounds (with high probability) under such constraints.

\paragraph{Fast rates in the batch setting:}
Another line of works investigated the problem of experts (or estimators) aggregation in the batch setting with stochastic and i.i.d samples (i.e., each expert's predictions are assumed to follow an independent and identical distribution, see \citealp{tsybakov2003optimal}). There are two distinct phases: a first step where the learner has access to training data points, then a prediction step where she outputs a combination of experts. The output in this setting is compared against the best expert. A non-exhaustive list of works considering this problem includes those of \cite{audibert2007progressive}, \cite{lecue2009aggregation}, and \cite{saad2021fast}, where the emphasis was put on obtaining $\mathcal{O}(1/T)$ ``fast'' rates for excess risk with high probability under some convexity assumptions on the loss function. However, these algorithms are not translatable to the adversarial setting since some of the previous strategies rely on the early elimination of sub-optimal experts. \cite{saad2021fast} presented a budgeted setting where the learner is constrained to see at most $m$ experts forecasts per data point and can predict using $p$ experts. This paper is an extension of their framework in the adversarial setting with a cumulative regret. 

\paragraph{Online Convex Optimization with bandit feedback:}
A different objective is considered in the online convex optimization framework, where the losses are compared against the best convex combination of the experts. This problem was studied by \cite{agarwal2010optimal} and \cite{shamir2017optimal} under limited feedback. More precisely, the learner can query the value of the loss function in two points from the convex envelope of the compact set over which the optimization is performed. In such a setting, it was shown that for Lipschitz and strongly-convex loss functions, it is possible to achieve an expected regret bounded by $\mathcal{O}\left(d^{2}\log(T)\right)$, where $d$ is the dimension of the linear span of experts (which plays a similar role to $K$ in our setting). Observe that online convex optimization algorithms (eg. as considered in the cited references) cannot be applied in our setting, where the player is not allowed to play (or observe) an arbitrary point in the convex envelope of the experts, but rather convex combinations with support on $p$ (or $m$) experts. On the other hand, the goal aimed at is different as well, since we want to minimize the regret with respect to the best expert, not with respect to the best convex combination of experts (which would not be an attainable goal under the considered play restrictions).

\paragraph{Why aim at high probability bounds instead of expectation bounds?} Consider an algorithm with internal randomization. From a practical point of view, bounds on its expected regret do not necessarily translate into a similar guarantee with high probability. In many applications, such as finance, controlling the fluctuations of risk is very important. From a mathematical point of view, the ``phenomenon" of negative regrets occurs when the player has a chance of outperforming the benchmark (such as the best-fixed expert in hindsight) for some rounds. In this case, an algorithm may have optimal expected regret but sub-optimal deviations. A manifestation of this problem is for the EXP3 algorithm in multi-armed bandit setting ($p=m=1$ in Protocol~\ref{algo:gp}), which has a worst case regret of $\sqrt{KT}$ in expectation, but the random regret can be linear $\Omega(T)$ with constant probability (see the exercises of Chapter 11 of \citealp{lattimore2020bandit}).

\section{Main results: Algorithm with upper bounds in expectation}\label{sec:Main1}

In this section, we introduce a new algorithm with constant bounds on the expected regret, for the setting:  $p= 2$ and $m\ge 3$. The aim of this section is to present some central intuitions, which are complemented in the next section to achieve stronger guarantees. To ease notation, we denote for each $i\in \intr{K}$ and $t\in \intr{T}$: $\ell_{i,t} := \ell\left(F_{i,t}, y_t\right)$. 

The high-level idea of Algorithm~\ref{algo:0} is common in the literature. It consists in constructing unbiased estimates of unseen losses, %such that $\mathbb{E}[\hat{\ell}_{i,t}| \mathcal{F}_{t-1}] = \ell_{i,t}$,
which are fed to the classical exponential weighting (EW) scheme over the experts. The first novelty introduced here is that the estimates are centered in a ``data-dependent" way, whose goal is to reduce variance. This variance control is essential in our analysis (see sketch of the proof below) in order to have constant regrets.

Let us denote $\hat{p}_t$ the probability distribution derived by the EW principle
using estimated cumulated losses $\hat{L}_{i,t}$ over the set of experts at round $t$.
The second novelty consists in sampling just two experts $I_t$ and $J_t$, independently at random following $\hat{p}_t$, and $m-2$ additional experts uniformly at random for exploration. Then, we play the mid-point of the predictions of $I_t$ and $J_t$ (i.e., predict we predict $\frac{1}{2}F_{I_t,t}+\frac{1}{2}F_{ J_t,t}$). 

%Loss functions satisfying Assumption~\ref{assump} satisfy inequality \eqref{eq:prop}. 
The main idea for getting a constant regret bound is to compensate the variance term introduced by the estimates $(\hat{\ell}_{i,t})$ by the negative second order term in inequality~\eqref{eq:prop} satisfied by the loss. 
The following theorem presents a constant bound on the expected regret, with a sketch of the proof.

\begin{algorithm}
	%\centering
	\caption{Prediction with limited advice $(p=2, m \ge 3)$}
	\begin{algorithmic}\label{algo:0}
		\STATE \textbf{Input Parameters:} $\lambda$, $m$.
		\STATE \textbf{Initialize:} $\hat{L}_{i,0} = 0$ for all $i \in \intr{K}$.
		\FOR{\textbf{each} round $t=1,2,\dots$}
		\STATE Let 
		\begin{equation*}
		\hat{p}_{i,t} = \frac{\exp\left( -\lambda \hat{L}_{i,t-1} \right)}{\sum_{j} \exp\left( -\lambda \hat{L}_{j, t-1} \right)}.
		\end{equation*}
		\STATE Draw $I_t$ and $J_t$ according to $\hat{p}_t$ independently.
		\STATE Play: $\frac{1}{2}F_{I_t,t}+\frac{1}{2}F_{ J_t,t}$, and incur its loss.
		\STATE Sample $m-2$ experts uniformly at random without replacement from $\intr{K}$. Denote $\cU_t$ this set of experts. 
		\STATE Query $C_t = \cU_t \cup \{I_t, J_t\}$.
		%\STATE \textbf{Decoupled exploration-exploitation:} Sample $m-1$ experts uniformly at random without replacement from $\intr{K}$. Denote $\tilde{\mathcal{O}}_t$ this set of experts and let $\tilde{m} = \lvert \tilde{\mathcal{O}}_t \rvert$. Query $C_t = \{I_t \} \cup \tilde{\mathcal{O}}_t$.
		\FOR{ $i \in \intr{K}$}
		%\STATE Let $\hat{q}_{i,t} =  \hat{p}_{i,t} + \hat{p}_{i,t} (1 - \hat{p}_{i,t}) +(1 - \hat{p}_{i,t})^2~\frac{m-2}{K-2}$.
		\STATE Let
		\begin{equation*}
		\hat{\ell}_{i,t} = \frac{K}{m-2}\mathds{1}\left(i \in \cU_t\right)~ \ell_{i,t} + \left( 1 - \frac{K}{m-2}\mathds{1}\left(i \in \cU_t\right)\right)~ \ell_{I_t,t}.
		\end{equation*}
		\STATE Update $\hat{L}_{i,t} = \hat{L}_{i, t-1} + \hat{\ell}_{i,t}$.  
		\ENDFOR  
		\ENDFOR
	\end{algorithmic}
\end{algorithm}

\noindent Define the following constant
\begin{equation}\label{eq:c_cst}
\bar{\lambda} := \min\left\lbrace\frac{4\log\left(1 + \frac{\eta^2 B^2}{2}\right)}{\eta B^2},~\frac{1}{B} \right\rbrace.
\end{equation}

\begin{theorem}\label{thm:0}
	Suppose Assumption~\ref{assump} holds. For any input parameter: $\lambda \in \left(0, \frac{m-2}{4K}\bar{\lambda}\right)$, where $\bar{\lambda}$ is defined in \eqref{eq:c_cst}, the expected regret of Algorithm~\ref{algo:0} satisfies:
	\begin{equation*}
	\mathbb{E}\left[ \mathcal{R}_T \right] \le \frac{\log(K)}{\lambda},
	\end{equation*}
	where the expectation is with respect to the learner's own randomization.
\end{theorem}
\begin{remark}
	Comparing this result with the guarantees of the classical exponential weights averaging (EWA) algorithm, one can notice that in the full information feedback setting ($m=K$), our guarantee is of the same order, up to a numerical constant, as the constant regret bound for EWA for exp-concave losses. The advantage of our procedure is that it necessitates sampling only two experts from the EW distribution instead of full averaging. In the partial feedback case ($m<K$),  Algorithm~\ref{algo:0} guarantees a regret of order $\cO(K\log(K)/m)$, as one would expect, the factor $K/m$ reflects the proportion of the information available to the learner. The last bound is tight, up to a logarithmic factor in $K$ (see Theorem~\ref{thm:2}).
	%Notice that the introduced estimator $\hat{\ell}_{i,t}$ of unseen losses are built using a centering technique, whose goal is to reduce variance in a data-dependent way. This variance control is essential in our analysis (see sketch of the proof below) in order to have constant regrets. 
\end{remark}

\begin{sproof}
	Let $\left( \mathcal{F}_t \right)$ denote the natural filtration associated to the process
	of available information, $\left(S_t, C_t, (\ell_{i,t})_{t \in C_t}\right)$, and denote $\mathbb{P}_{t-1}$ resp. $\mathbb{E}_{t-1}$ the conditional
	probability resp. expectation with respect to $\mathcal{F}_{t-1}$ (``past observations''). The loss functions $\ell_t$ satisfy Assumption~\ref{assump}. Therefore, using Lemma~\ref{lem:assumption1}, 
	the expected cumulative loss of Algorithm~\ref{algo:0} is given by
	\begin{align}
	\sum_{t=1}^{T} \mathbb{E}\left[\ell_t\left(\frac{F_{I_t,t}+F_{J_t,t}}{2}\right)\right] &\le \sum_{t=1}^{T} \mathbb{E}\left[\frac{1}{2}\ell_{I_t,t} + \frac{1}{2}\ell_{J_t,t} - \frac{\bar{\lambda}}{2} (\ell_{I_t,t}-\ell_{J_t,t})^2\right] \notag\\
	& =\underbrace{\sum_{t=1}^{T}\sum_{i=1}^{K} \mathbb{E}\left[\hat{p}_{i,t}~\ell_{i,t}\right]}_{\text{Term 1}} - \underbrace{\frac{\bar{\lambda}}{2} \sum_{t=1}^{T}\sum_{i,j=1}^{K}\mathbb{E}\left[\hat{p}_{i,t}\hat{p}_{j,t} \left(\ell_{i,t} - \ell_{j,t}\right)^2\right]}_{\text{Term 2}}.\label{eq:bound_1}
	\end{align}
	
	%		\noindent $\ell_t$ satisfies Assumptions~\ref{assump}. Therefore, using Lemma~\ref{lem:assumption1}, it satisfies inequality \eqref{eq:prop}. 
	%		\begin{equation}\label{eq:bound_1}
	%		\sum_{t=1}^{T} \mathbb{E}\left[\ell_t\left(\frac{F_{I_t,t}+F_{J_t,t}}{2}\right)\right] = \underbrace{\sum_{t=1}^{T}\sum_{i=1}^{K} \mathbb{E}\left[\hat{p}_{i,t}~\ell_{i,t}\right]}_{\text{Term 1}} - \underbrace{\frac{\bar{\lambda}}{2} \sum_{t=1}^{T}\sum_{i,j=1}^{K}\mathbb{E}\left[\hat{p}_{i,t}\hat{p}_{j,t} \left(\ell_{i,t} - \ell_{j,t}\right)^2\right]}_{\text{Term 2}}.
	%		\end{equation} 
	
	\noindent Observe that by construction of Algorithm~\ref{algo:0}, the elements in $\cU_t$ were sampled uniformly at random without replacement  from $\intr{K}$. Moreover, $\cU_t$ is independent of $I_t$. Therefore, $\hat{\ell}_{i,t}$ is an unbiased estimator of $\ell_{i,t}$
	conditionally to the available information: $\mathbb{E}_{t-1}[ \hat{\ell}_{i,t} ] = \ell_{i,t}$.
	%        since
	%	\begin{align*}
	%	\mathbb{E}_{t-1}\left[ \hat{\ell}_{i,t} \right] &= \mathbb{E}_{t-1}\left[ \frac{K}{m-2}\mathds{1}\left(i \in \cU_t\right)~ \ell_{i,t} + \left( 1 - \frac{K}{m-2}\mathds{1}\left(i \in \cU_t\right)\right)~ \ell_{I_t,t}\right]\\
	%	&=\frac{K}{m-2}\mathbb{P}_{t-1}\left[ i \in \cU_t \right] \ell_{i,t} + \paren{1 -  \frac{K}{m-2} \mathbb{P}_{t-1}\left[ i \in \cU_t\right]} \mathbb{E}_{t-1}\left[\ell_{I_t,t}\right]\\
	%	&= \ell_{i,t},
	%	\end{align*}
	%	where we used the independence of $\cU_t$ and $I_t$ in the second line,
	%        and $\mathbb{P}_{t-1}\left[ i \in \cU_t\right] = (m-2)/K$ in the last line.

	Using the tower rule, Term 1 therefore writes $\sum_{t}\sum_{i} \mathbb{E}[\hat{p}_{i,t} \hat{\ell}_{i,t}]$. Next, we use Lemma~\ref{lem:log_tel} in the Appendix (by cancellation of consecutive logarithmic terms) with $\mu_t = \sum_{i=1}^{K} \hat{p}_{i,t} \ell_{i,t}$ for each $t \in \intr{T}$. We have the following upper bound for Term 1 in \eqref{eq:bound_1}:
	\begin{equation}\label{eq:bound2}
	\sum_{t=1}^{T}\sum_{i=1}^{K} \mathbb{E}\left[\hat{p}_{i,t}~\hat{\ell}_{i,t}\right] \le \min_{i \in \intr{K}} \sum_{t=1}^{T} \mathbb{E}\left[\hat{\ell}_{i,t}\right] + \frac{\log(K)}{\lambda} + \lambda \sum_{t=1}^{T} \sum_{i=1}^{K} \mathbb{E}\left[\hat{p}_{i,t} \left(\hat{\ell}_{i,t} - \mu_t\right)^2\right].
	\end{equation} 
	
	\noindent We use the definition of $\hat{\ell}_{i,t}$ and the tower rule to upper bound the last term in \eqref{eq:bound_1}:
	\begin{align*}
	\mathbb{E}\left[ \sum_{i=1}^{K} \hat{p}_{i,t} \left(\hat{\ell}_{i,t} - \mu_t\right)^2\right] &\le \frac{2K}{m-2} \mathbb{E}\left[\sum_{i=1}^{K} \hat{p}_{i,t} \left(\ell_{i,t} - \mu_t\right)^2 \right] + \frac{2K}{m-2} \mathbb{E}\left[\left(\ell_{I_t,t} - \mu_t\right)^2\right]\\
	&= \frac{4K}{m-2} \mathbb{E}\left[\sum_{i=1}^{K} \hat{p}_{i,t} \left(\ell_{i,t} - \mu_t\right)^2\right].
	\end{align*}
	
	\noindent Finally, we combine \eqref{eq:bound_1}, \eqref{eq:bound2} and the bound above to obtain
	\begin{equation*}
	\mathbb{E}\left[\mathcal{R}_T\right] \le \frac{\log(K)}{\lambda} +\lambda \frac{4K}{m-2}  \mathbb{E}\left[\sum_{i=1}^{K} \hat{p}_{i,t} \left(\ell_{i,t} - \mu_t\right)^2\right] - \bar{\lambda} \sum_{t=1}^{T}\sum_{i,j=1}^{K}\mathbb{E}\left[\hat{p}_{i,t}\hat{p}_{j,t} \left(\ell_{i,t} - \ell_{j,t}\right)^2\right].
	\end{equation*}
	
	\noindent Recall that if $X$ and $Y$ are two independent and identically distributed variables, we have $\mathbb{E}[(X-Y)^2] = 2\Var(X)$. Applying this identity to Term 2 in \eqref{eq:bound_1}, we have
	
	\begin{equation*}
	\mathbb{E}\left[\mathcal{R}_T\right] \le \frac{\log(K)}{\lambda} +\left(\lambda \frac{4K}{m-2} - \frac{1}{B} \right)  \mathbb{E}\left[\sum_{i=1}^{K} \hat{p}_{i,t} \left(\ell_{i,t} - \mu_t\right)^2\right].
	\end{equation*}
	We conclude using $\lambda < \frac{m-2}{4K} \bar{\lambda}$. 
\end{sproof}

\section{Main results: Algorithms with high probability upper bounds}\label{sec:Main2}

In this section, we present new algorithms with guarantees that hold with high probability with respect to the player's own randomization. As discussed in Section~\ref{sec:related_w}, high probability guarantees are important to assess any algorithm's goodness due to potential exposure to negative regrets phenomena and
thus the possibility of deviations having larger order than the expectation.

We introduce sampling strategies for three different settings:  $p= 2$ and $m\ge 3$, $(p=2, m=2, \text{IC}=\text{False})$ and $(p=2, m=2, \text{IC}=\text{True})$, presented in Algorithms~\ref{algo:1}
and~\ref{algo:2}; Algorithm~\ref{algo:1} is common to the first two settings. To ease notations, we denote for each $i\in \intr{K}$ and $t\in \intr{T}$: $\ell_{i,t} := \ell\left(F_{i,t}, y_t\right)$.
%For simplicity, we will discuss here only the settings $(p= 2, m\ge3)$ and $(p=2, m=2, \text{IC}=\text{True})$.

In Algorithms~\ref{algo:1} and~\ref{algo:2}, we build on the idea presented in Algorithm~\ref{algo:0} and construct estimates of unseen losses, which are fed into an EW scheme from which experts are sampled.
%The main difference with
%Algorithm~\ref{algo:0} is that we additionally introduce a {\em negative} (or ``optimistic'') bias on the estimated losses, which takes into account an estimated variance. This can be conceptually compared
%to the uniform confidence bound (UCB) algorithm in the standard bandit setting, which will select
%``optimistically'' arms which have the highest potential reward given past information (here, loss is a negative reward).
Let $\hat{p}_t$ denotes the resulting estimated EW distribution.
The main differences between the algorithms below and Algorithm~\ref{algo:0} are (a) the constructed loss estimates and (b) the sampling strategy when $m=2$ and $\text{IC}=\text{True}$.

{\bf Modified loss estimates:}
%We introduce unbiased estimators for unseen losses $\mathbb{E}[\hat{\ell}_{i,t}| \mathcal{F}_{t-1}] = \ell_{i,t}$. These estimates are centered in a ``data-dependent" way at $\ell_{A_t,t}$ (where $A_t \sim \hat{p}_t$).
We start with the same unbiased loss estimates, with data-dependent centering, from
Algorithm~2, but additionally introduce a {\em negative} (or ``optimistic'') bias on the estimated losses, which takes into account an estimated variance. This can be conceptually compared
to the uniform confidence bound (UCB) algorithm in the standard stochastic bandit setting, which will select
``optimistically'' arms which have the highest potential reward given past information (here, loss is a negative reward). In this sense, this term tends to encourage diversity in expert sampling (i.e.
encourage sampling experts with a possibly higher estimated loss but also larger variance than the
best estimated experts so far). 
This is used in both Algorithms~\ref{algo:1} and~\ref{algo:2}.
%a second order term $\hat{V}_{i,t}$ to these estimates to encourage playing experts that are distant. The resulting biased estimate is then fed to the classical exponential weights. 

In the case $m\ge 3$ or $(m=2, \text{IC}=\text{False})$, there is still at least one free observation
left for exploration decoupled from exploitation. In these settings, Algorithm~\ref{algo:1} uses
the same sampling scheme as Algorithm~\ref{algo:0}, % the playing scheme consists of
namely sampling independently at random two experts following $\hat{p}_t$ and playing the central point of the sampled predictions. The remaining ``pure exploration'' observations are sampled uniformly at random, with replacement.
%, as in Algorithm~\ref{algo:0}.

{\bf Modified sampling scheme:} the case $(m=2, \text{IC}=\text{True})$ is more difficult since there
is no ``free exploration'' observation possible. This is the counterpart of the 
exploration/exploitation tradeoff of the standard bandit setting, in
the framework where we aim at constant regrets (so that playing combinations of at least two arms is necessary, see next section). Taking inspiration from the standard bandit setting literature ($p=m=1$),
introducing a small uniform exploration component appears necessary for the sampling strategy for algorithms achieving optimal high probability guarantees (\citealp{audibert2010regret}, \citealp{auer2002nonstochastic}, \citealp{beygelzimer2011contextual}, \citealp{bubeck2012regret}). For example, EXP3.P mixes the EW sampling rule with a uniform distribution over the arms. On the other hand, EXP-IX \citep{neu2015explore} incorporates the exploration component implicitly through a biased estimate of the losses. However, this uniform exploration costs $\mathcal{O}(\sqrt{KT})$ on the cumulative regret. Hence, aiming at constant regret necessitates a more subtle sampling rule. 

We introduce a two-step sampling strategy. The first expert, denoted $A_t$, is sampled following $\hat{p}_t$. The second expert, denoted $B_t$, is sampled uniformly at random (possibly $B_t$ and $A_t$ are identical). The predictions of $(A_t,B_t)$ are observed after making a prediction. For the playing strategy, we sample two experts independently (conditionally to $A_t$ and $B_t$) at random, following the restriction of the law $\hat{p}_t$ on  $\{A_t,B_t\}$, and we play the central point of the two sampled experts. Therefore, depending on the outcome of the second step, the algorithm's prediction can be either one of the two pre-selected experts or the central point of the two experts. This strategy ensures the necessary uniform exploring component needed in the adversarial problems.

The possibility of having constant regrets guarantees is due to Property \eqref{eq:prop}, satisfied for the loss functions $\ell$ under Assumption~\ref{assump}: Lemma~\ref{lem:assumption1} suggests that when predicting the central point of two experts, the learner benefits from the distance between the played predictions. This remark is exploited in constructing of the distribution $\hat{p}_t$.
%$(\hat{p}_{i,t})_{i,t}$. 

To summarize, the playing strategy relies on three essential ideas: the (conditional for $m=2$) independence of the played experts, the centering scheme for the losses estimates, and the second order term to diversify the played arms.

\begin{algorithm}[!ht]
	% \centering
	\caption{\label{algo:1}  ($p = 2$, $m \ge 3$) or ($p = 2$, $m = 2$, $\text{IC}=\text{False}$) }
	\begin{algorithmic}
		\STATE \textbf{Input Parameters:} $\lambda,m$.
		\STATE \textbf{Initialize:} $\hat{L}_{i,0} = 0, \hat{V}_{i,0} = 0$ for all $i \in \intr{K}$.
		\STATE Let $\tilde{m} = \max\{m-2,1 \}$.  
		\FOR{\textbf{each} round $t=1,2,\dots$}
		\STATE Let 
		\begin{equation}\label{eq:defphat}
		\hat{p}_{i,t} = \frac{\exp\left( -\lambda \hat{L}_{i,t-1} + \lambda^2 \hat{V}_{i,t-1} \right)}{\sum_{j=1}^K \exp\left( -\lambda \hat{L}_{j, t-1} + \lambda^2 \hat{V}_{j,t-1} \right)}.
		\end{equation}
		\STATE Sample $I_t$ and $J_t$ according to $\hat{p}_t$ from $\intr{K}$ independently.
		\STATE \textit{Play:} $\frac{1}{2}F_{I_t,t}+\frac{1}{2}F_{ J_t,t}$, and incur its loss.
		\STATE Sample $\tilde{m}$ experts without replacement, independently and uniformly at random from $\intr{K}$. Denote $\cU_t$ this set of experts. 
		\IF{$m\ge3$}
		\STATE Let $C_t = \{I_t,J_t\}\cup \cU_t$.
		\ELSIF{$m=2$}
		\STATE Let $C_t = \{I_t\}\cup \cU_t$.
		\ENDIF
		\STATE \textit{Observe:} $\ell_{i,t}$ for $i \in C_t$.
		\FOR{ $i \in \intr{K}$}
		\STATE Let
		\begin{align}\label{eq:deflhat}
		\hat{\ell}_{i,t} &= \frac{K}{\tilde{m}}\mathds{1}\left(i \in \mathcal{U}_t\right)~ \ell_{i,t} + \left( 1 - \frac{K}{\tilde{m}}\mathds{1}\left(i \in \mathcal{U}_t\right)\right)~ \ell_{I_t,t}\\
		\hat{v}_{i,t} &= \left(\hat{\ell}_{i,t} - \ell_{I_t,t}\right)^2 \label{eq:defvhat}
		\end{align}
		\STATE Update $\hat{L}_{i,t} = \hat{L}_{i, t-1} + \hat{\ell}_{i,t}$ and $\hat{V}_{i,t} = \hat{V}_{i, t-1} + \hat{v}_{i,t}$.  
		\ENDFOR  
		\ENDFOR
	\end{algorithmic}
\end{algorithm}

\begin{algorithm}[!ht]
	% \centering
	\caption{\label{algo:2}  ($p = 2$, $m = 2$, $\text{IC}=\text{True}$) }
	\begin{algorithmic}
		\STATE \textbf{Input Parameters:} $\lambda$.
		\STATE Initialize: $\hat{L}_{i,0} = 0$ for all $i \in \intr{K}$.
		\FOR{\textbf{each} round $t=1,2,\dots$}
		\STATE Let 
		\begin{equation*}
		\hat{p}_{i,t} = \frac{\exp\left( -\lambda \hat{L}_{i,t-1} + \lambda^2 \hat{V}_{i,t-1} \right)}{\sum_{j =1}^K \exp\left( -\lambda \hat{L}_{j, t-1} + \lambda^2 \hat{V}_{j,t-1} \right)}.
		\end{equation*}
		\STATE Sample one expert from $\intr{K}$, denoted $A_t$, according to $\hat{p}_t$, and one expert from $\intr{K}$, denoted $B_t$, independently and uniformly at random. Let $C_t = \{A_t, B_t\}$.
		\FOR{$i \in C_t$} 
		\STATE Let 
		\begin{equation*}
		\hat{q}_{i,t} = \frac{\exp\left( -\lambda \hat{L}_{i,t-1} + \lambda^2 \hat{V}_{i,t-1} \right)}{\sum_{j \in C_t}\exp\left( -\lambda \hat{L}_{j, t-1} + \lambda^2 \hat{V}_{j,t-1} \right) }.
		\end{equation*}
		\STATE Draw $I_t$ from $C_t$ according to $\hat{q}_t$.
		\STATE \textit{Draw} $J_t$ from $C_t$ according to $\hat{q}_t$ independently from $I_t$.
		%\IF{$p>2$} 
		%\STATE\textit{Draw} $J_t$ from $\intr{K}$ according to $\hat{p}_t$ independently from $I_t$.
		%\ENDIF
		\STATE \textit{Play:} $\frac{1}{2}F_{I_t,t}+\frac{1}{2}F_{ J_t,t}$, and incur its loss.
		\STATE \textit{Observe:} $\ell_{i,t}$ for $i \in C_t$.
		\ENDFOR
		\FOR{ $i \in \intr{K}$}
		\STATE Let
		\begin{align*}
		\hat{\ell}_{i,t} &= K~\mathds{1}\left(B_t=i\right)~ \ell_{i,t} + \left( 1 - K~\mathds{1}\left(B_t = i\right)\right)~ \ell_{A_t,t}\\
		\hat{v}_{i,t} &= \left(\hat{\ell}_{i,t} - \ell_{A_t,t}\right)^2
		\end{align*}
		\STATE Update: $\hat{L}_{i,t} = \hat{L}_{i, t-1} + \hat{\ell}_{i,t}$ and $\hat{V}_{i,t} = \hat{V}_{i, t-1} + \hat{v}_{i,t}$.  
		\ENDFOR  
		\ENDFOR
	\end{algorithmic}
\end{algorithm}

\begin{remark}
	\begin{itemize}
		\item The proposed algorithm can be implemented in an efficient way, so that after a one-time computational
		cost of $\cO(K)$ for initialization, the computational cost of each round, including suitably keeping track
		of the distribution $\hat{p}_t$ and sampling from it, is $\cO(m \log K)$ (see Appendix~\ref{app:algocomplex} for details). Therefore, the computational complexity also depends mildly on the number of experts $K$.
		% Complexity of sampling from $\hat{p}_t$ is $O(\log K)$: this holds for sampling from a distribution on a finite set of cardinality $K$ if  states are organized as leaves $\partial T $ of a balanced binary tree $T$ storing
		% quantities of the form $S_{t_0} = C.\sum_{t \in \partial T_{t_0}}\hat{p}_t$, where $T_{t_0}$ is the subtree rooted
		% at node $t_0$. Namely, put $Z=S_{\mathrm{root}}. U = CU$ where $U \sim \mathrm{Unif}[0,1]$, then
		% go left or right down the tree recursively according to $Z< S_{\mathrm{left}}$ or not,
		% and update $Z \leftarrow Z - S_{\mathrm{left}}$ if going right.
		\item  Since our analysis suggests that we can restrict possible plays to mid-points of just two experts, one could argue that the coupled setting $(p = m = 2 , \text{IC=True})$ looks quite similar to learning with expert advice with bandit
		feedback, where the possible arms would be the $K^2$ “bi-experts” that are mid-points of original experts $(i, j)$.
		One could therefore think of a more direct approach: simply applying a bandit-type strategy, say EXP3.P or EXP3-IX (\citealp{auer2002nonstochastic} and \citealp{neu2015explore}, respectively) to these $K^2$ ``arms''.   However,
		existing generic results only guarantee a “slow” $\mathcal{O}(\sqrt{T})$ regret with respect to the best “bi-expert”, and
		this cannot be compensated in general by exp-concavity, as the best “bi-expert” may not be much
		better than the best expert (if the experts are “correlated”: see proof of lower bounds in Theorem~\ref{thm:2} and~\ref{thm:4}). Furthermore, in the playing strategy of EXP3.P and EXP3-IX, each pair of experts is played $\Omega(\sqrt{K^2T})$ times, due the uniform exploration component of their sampling schemes. This will lead regrets scaling with $\sqrt{T}$. 
	\end{itemize}

\end{remark}

%	\begin{proposition}\label{prop:sample}
%		Consider Algorithm~\ref{algo:1} in the case $\bm{p=m=2}$. For each $t\in \intr{K}$ and $i\in \intr{K}$, we have
%		\[
%		\mathbb{P}\left(I_t = i | \mathcal{F}_{t-1}\right) = \hat{p}_{i,t}.
%		\]
%		Moreover
%		\[
%		\mathbb{P}\left(I_t = i, J_t = j | \mathcal{F}_{t-1}\right) \ge \frac{1}{K} \hat{p}_{i,t}\hat{p}_{j,t}.
%		\]
%	\end{proposition}

%Below we present the main theorems of this work. 

\begin{theorem}\label{thm:1}
	Suppose Assumption~\ref{assump} holds. 
	
	Consider the case $ (m\ge 3 \text{ and } p = 2) \text{ or } (m= 2 \text{ and } p = 2 \text{ and } \text{IC}=\text{False})$. For any input parameter $\lambda \in \left(0, \frac{m-1}{128K} \bar{\lambda}\right)$, where $\bar{\lambda}$ is defined in \eqref{eq:c_cst}, the regret of Algorithm~\ref{algo:1} satisfies with probability at least $1 - 8 \delta$, with respect to the player's own randomization
	\begin{equation*}
	\mathcal{R}_T \le c~\frac{1}{\lambda}  \log\left(\frac{\bar{\lambda}K}{\lambda  \delta}\right),
	\end{equation*}
	where $c$ is a numerical constant.
\end{theorem}

\begin{theorem}\label{thm:1bis}
	Suppose Assumption~\ref{assump} holds. 
	
	Consider the case $p=m=2$ and $\text{IC}=\text{True}$. For any input parameter $\lambda \in \left(0, \frac{\bar{\lambda}}{352K^2}\right)$, where $\bar{\lambda}$ is defined in \eqref{eq:c_cst}, the regret of Algorithm~\ref{algo:2} satisfies with probability at least $1 - 8 \delta$, with respect to the player's own randomization
	\begin{equation*}
	\mathcal{R}_T  \le c~\left(\frac{1}{\lambda} + \frac{K}{\bar{\lambda}}\right) \log\left(\frac{\bar{\lambda}K}{\lambda  \delta}\right),
	\end{equation*}
	where $c$ is a numerical constant.
\end{theorem}

\paragraph{Discussion} Notice that prior knowledge on the confidence level $\delta$ is not required by Algorithms~\ref{algo:1} and~\ref{algo:2}. The presented bounds in theorems above are valid for any $\delta \in (0,1)$. Observe that taking $\lambda$ close to $m/(128K)~ \bar{\lambda}$ leads to a bound of the order $\cO(K\log(K\delta^{-1})/m)$ in Theorem~\ref{thm:1}, which is minimax optimal up to a $\log(K)$ factor (Theorem~\ref{thm:2}). Taking $\lambda$ close to $1/(352 K^2)~\bar{\lambda}$, leads to a bound of the order $\cO(K^2\log(K\delta^{-1}))$ in the special setting $p=m=2$ with $\text{IC} = \text{True}$. This bound presents a gap of factor $K$ with the lower bound presented in Theorem~\ref{thm:2}. We emphasize that in the last setting, the player chooses two experts to combine their predictions and observes only the feedback of these two experts. Hence, unlike the setting considered in Theorem~\ref{thm:1}, the player is deprived of additional 'freely chosen' experts to explore their losses. This constraint necessitates a more careful playing strategy, presented in Algorithm~\ref{algo:2}.   

\section{Lower bounds}\label{sec:lower}
In this section, we provide lower bounds matching the upper bounds in Theorem~\ref{thm:1}, up to a logarithmic factor in $K$ (except for the case $p=m=2$, where we have a gap of factor $K$). The techniques of the proof are similar to the ones presented by \cite{auer1995gambling}. The main difference comes from the construction of the experts' distributions.
%\textcolor{red}{TODO:} get rid of $K\ge 5$ in the proof of some theorems.
\begin{theorem}\label{thm:2}
	Let $\ell$ be the squared loss:  $\ell(x,y)=(x-y)^2$ on $\cX=\cY=[0,1]$.
	Consider the game protocol presented in Algorithm~\ref{algo:gp} with $m\ge 2 \text{ and } p \ge 2 \text{ and } \text{IC} \in \{\text{False},\text{True} \}$. The expected regret satisfies:
	\begin{equation*}
	\inf \sup \mathbb{E}\left[ \mathcal{R}_T \right] \ge c~\frac{K}{m},
	\end{equation*}
	where $c$ is a numerical constant, the infinimum is over all playing strategies and the supremum is over all individual sequences.
\end{theorem}
\begin{remark}
	The lower bound presented in Theorem~\ref{thm:2} is valid for any $p \le K$. Algorithms~\ref{algo:1} and~\ref{algo:2} match it (up
	to a log factor in K) using only $p = 2$, suggesting that no significant improvements can be obtained
	if we are allowed to predict using more than two experts.
\end{remark}

Theorem below is of theoretical interest, it shows that if only one feedback is received per round, then constant regrets are not achievable.
\begin{theorem}\label{thm:4}
	Let $\ell$ be the squared loss:  $\ell(x,y)=(x-y)^2$ on $\cX=\cY=[0,1]$.
	Consider the game protocol presented in Algorithm~\ref{algo:gp} with $m=1 \text{ and } p\in \intr{K} \text{ and } \text{IC} \in \{\text{False},\text{True}\}$, we have
	\begin{equation*}
	\inf \sup \mathbb{E}\left[ \mathcal{R}_T \right] \ge c~\sqrt{KT},
	\end{equation*}
	where $c$ is a numerical constant, the infinimum is over all playing strategies and the supremum is over all individual sequences.
\end{theorem}

For the sake of completeness, we state the following lower bound from \cite{seldin2014prediction}.

\begin{theorem}[Direct consequence of \citealp{seldin2014prediction}]\label{thm:5}
	Let $\ell$ be the squared loss:  $\ell(x,y)=(x-y)^2$ on $\cX=\cY=[0,1]$.
	Consider the game protocol presented in Algorithm~\ref{algo:gp} with $p=1 \text{ and } m\in \intr{K} \text{ and } \text{IC} \in \{\text{False},\text{True}\}$, we have
	\begin{equation*}
	\inf \sup \mathbb{E}\left[ \mathcal{R}_T \right] \ge c~\sqrt{\frac{K}{m}T},
	\end{equation*}
	where $c$ is a numerical constant, the infinimum is over all playing strategies and the supremum is over all individual sequences.
\end{theorem}

\section{Discussion and open questions}

\begin{itemize}
	\item In the setting $p=m=2$ with coupled exploration-exploitation ($\text{IC}=\text{True}$), Algorithm~\ref{algo:2} presents a strategy with a bound of order $\cO(K^2\log(K\delta^{-1}))$, while the lower bound presented in Theorem~\ref{thm:2} is of order $\cO(K)$. It would be of interest to close this gap.   
	\item Previous works on achieving constant regret under a full observation model only assumed exp-concavity of
	the loss (see e.g. \citealp{cesa2006prediction}, Chap. 3). In the limited observation setting, we additionally assume that the loss function is bounded by a constant $B$ known to the player. It would be of interest to determine if this condition is necessary. We note, however that loss boundedness is an important ingredient in applying Bernstein-type inequalities for bounds in high
	probability.    
	\item In the stochastic (i.i.d. experts and target variables) setting, a variation of the expert elimination
	strategy proposed by \cite{saad2021fast} (suitably adapted to tackle cumulative regret) can be shown to have
	fast rates for regret in an instance-free setting, as well as suitable instance-dependent performance bounds (i.e.,
	the bound depends on the average performance of experts and their correlation, eliminating clearly sub-optimal experts earlier). This a fairly different strategy from the exponential weighting variations
	proposed here.  In the bandit setting, \cite{seldin2014one} have proposed a strategy that reaches almost optimal bounds
	both in the stochastic
	and the adversarial settings. It would be interesting to investigate whether such an omnibus strategy exists.
	\item We have shown that $p=2$ is sufficient to get constant regret with respect to the best expert,
	using a strong convexity-type assumption on the loss. For $p=K$, for an exp-concave loss there
	exist strategies having constant regret with respect to the best convex combination of experts (e.g. \citealp{cesa2006prediction}, Theorem. 3.3), albeit with a $O(K)$ scaling of the regret.
	It would be interesting to study if ``intermediate'' situations exist, for example if it is possible
	to have constant regret with respect to $k$-combinations of experts using only $p=\cO(k)$ expert
	predictions.
\end{itemize}

\bibliographystyle{plainnat}
\bibliography{bib_data_base}

\vfill

\pagebreak

\appendix

{\bf \Large Appendix: detailed proofs}

\section{Notation}\label{sec:not}

The following notation pertains to all the considered algorithms, where $t$ is a given training round and $T$ is the game horizon:
\begin{itemize}
	\item For any $x >0$, let $\log_2^{+}(x) = \max\{0, \log_2(x)\}$.
	\item Let $\mathcal{R}_T$ denote the cumulative random regret of the player over $T$ rounds.
	\item Let $S_t$ denote the set of combined experts to make a prediction at round $t$.
	\item Let $C_t$ denote the set of observed experts after making the prediction at round $t$.
	\item For each $i \in S_t$, let $\alpha_{i,t}$ denote the weight of expert $i$ in the convex combination played in round $t$.
	\item Let $\left(\mathcal{F}_t\right)_{t}$ denote the natural filtration associated with the process $\left(S_t, C_{t}, (\ell_{i,t})_{i \in C_{t}}\right)_t$.
	\item Denote the conditional expectation with respect to $\mathcal{F}_t$ by $\mathbb{E}_{t}[.] = \mathbb{E}\left[.|\mathcal{F}_t\right]$.
	\item For each expert $i \in \intr{K}$, let $N_i$ denote the number of times the prediction of expert $i$ was observed during the game (over $T$ rounds).
	\item For each expert $i \in \intr{K}$, let $M_i$ denote the number of times the prediction of expert $i$ was used for prediction during the game (over $T$ rounds): $M_i := \abs{\{t \in \intr{T}:\quad i\in S_t\}}$.
	\item For each expert $i \in \intr{K}$, we define $\ell_{i,t} = \ell\left(F_{i,t}, y_t\right)$.
	\item Denote by $\ell_t: \mathcal{X} \to \mathbb{R}$ such that $\forall x \in \intr{X}: \ell_t(x) = \ell(x,y_t)$.
\end{itemize}

Notation associated to Algorithms~\ref{algo:1} and~\ref{algo:2}
\begin{itemize}
	\item Let $I_t$ and $J_t$ denote the experts used for prediction in round $t$.
	\item Let $\cU_t$ the set of experts queried for exploration
	(sampled uniformly without replacement from $\intr{K}$). In Algorithm~\ref{algo:2} let $\cU_t = \{B_t\}$.
	\item Let $\tilde{m} = \max\{1, m-2\}$.
	%\item Let $\hat{p}_{t}$ denote the distribution over $\intr{K}$, according to which $I_t$ and $J_t$ were (independently) sampled: $\forall i \in \intr{K}:\hat{p}_{i,t} = \mathbb{P}\left(I_t = i\right)$. 
\end{itemize}

\section{Some preliminary technical results}

The following device is standard (it is used for instance for proving Bennett's inequality).
\begin{lemma}\label{lem:bennett}
	Let $X$ be a random variable with finite variance, such that $X \le b$ almost surely for some $b>0$. For any $\lambda> 0$:
	\begin{equation*}
	\log\left(\mathbb{E}e^{\lambda X}\right) \le \lambda \mathbb{E}[X] + \frac{\phi(\lambda b)}{b^2} \mathbb{E}[X^2].
	\end{equation*}
	Where $\phi(x) = \exp(x) - 1 - x$.
\end{lemma}
\begin{proof}
	The function $x \mapsto x^{-2}\phi(x)$ is non-decreasing on $\mbr$. As a consequence, if $X\leq b$ a.s., for
	any $\lambda>0$ it holds
	$\exp(\lambda X) \leq \frac{\phi(\lambda b)}{b^2}  X^2 + 1 + \lambda X$, a.s. Taking the expectation, then applying the inequality
	$\log(1+t) \leq t$ yields the result.
\end{proof}

\begin{corollary}\label{cor:bennett}
	Let $X$ be a random variable with finite variance, such that $X \ge -b$ almost surely for $b>0$. For any $\lambda \in \left(0, \frac{1}{b}\right)$:
	\begin{equation*}
	\log\left(\mathbb{E}e^{-\lambda X}\right) \le -\lambda \mathbb{E}[X] + \lambda^2 \mathbb{E}[X^2].
	\end{equation*}
\end{corollary}
\begin{proof}
	This corollary is a direct consequence of applying Lemma~\ref{lem:bennett} to the variable $-X \le b$, then using the fact that $\forall x \le 1: \phi(x) \le x^2$.
\end{proof}

We now introduce some technical lemmas used in the proofs. Let us start by reminding the following standard result (see Theorem 1.1.4 \citealp{niculescu2006convex}).

\begin{lemma}\label{lem:conv}
	A continuous function $f: \mathcal{X} \to \mathbb{R}$, where $\cX$ is a convex set, is convex if and only if: for any $x,x' \in \mathcal{X}$:
	\[
	f\left(\frac{x+x'}{2}\right) \le \frac{1}{2} f(x)+ \frac{1}{2} f(x').
	\] 
\end{lemma}

Lemmas below give some bounds for some functions.

\begin{lemma}\label{lem:colg}
	\begin{itemize}
		\item We have for any $x \in \mathbb{R}$
		\[
		1 + \frac{x^2}{2} \le \cosh(x) \le \exp(x^2/2).
		\]
		\item Let $c>0$. We have for any $x \in [0,c]$
		\[
		\log(1+x) \ge \frac{\log(1+c)}{c}x. 
		\]
	\end{itemize}
\end{lemma}
\begin{proof}
	The first and third result is a direct consequence of Taylor's expansion. The second result follows simply by concavity of $x \to \log(1+x)$.
\end{proof}

\begin{lemma}\label{lem:pure_cal}
	We have for any $x,y > 0$
	\[
	\log_2^{+}\left(x\right) - \frac{x}{y} \le  \log_2^{+}\left(y\right).
	\]
\end{lemma}
\begin{proof}
	Let $x,y>0$, we have 
	\begin{align*}
	\log_2(y) &= \log_2(x) - \log_2\left(\frac{x}{y}\right)\\
	&\ge \log_2(x) - \frac{x}{y},
	\end{align*}
	where we used the fact that $\log_2(t) \le t$ for any $t>0$. To conclude we use the inequality
	\[
	(a)_+ - b \le (a-b)_{+},
	\]
	valid for any $a \in \mathbb{R}$ and $b>0$.
\end{proof}

\section{Proof of Lemma~\ref{lem:assumption1}}

Let $y \in \mathcal{Y}$. In this proof, we will denote $\ell(.)$ instead of $\ell(.,y)$ so as to ease notation.

\subsection{First claim}

By exp-concavity of $\ell$, we have for any $x,x' \in \mathcal{X}$
\begin{equation*}
\frac{1}{2} \exp\left\lbrace-\eta \ell(x) \right\rbrace + \frac{1}{2} \exp\left\lbrace -\eta \ell(x') \right\rbrace\le \exp\left\lbrace-\eta  \ell\left(\frac{x+x'}{2}\right)\right\rbrace.
\end{equation*}
Multiplying both sides by $\exp\left\lbrace\frac{1}{2} \eta\ell(x) + \frac{1}{2} \eta\ell(x') \right\rbrace$, we have
\[
1 + \frac{\eta^2 \left(\ell(x) - \ell(x')\right)^2}{2}\le \exp\left\lbrace\frac{\eta}{2} \ell(x) + \frac{\eta}{2} \ell(x')-\eta  \ell\left(\frac{x+x'}{2}\right)\right\rbrace,	
\]
where we used the first result of Lemma~\ref{lem:colg} to lower bound the left hand side. 

\noindent Introducing the logarithm and using the second result of Lemma~\ref{lem:colg}, we obtain

\[
\frac{2 \log\left(1 + \frac{\eta^2 B^2}{2}\right)}{\eta^2 B^2} \eta^2 \left(\ell(x) - \ell(x')\right)^2 \le \frac{\eta}{2} \ell(x) + \frac{\eta}{2} \ell(x')-\eta  \ell\left(\frac{x+x'}{2}\right).
\]
We conclude that
\[
\ell\left(\frac{x+x'}{2}\right) \le \frac{1}{2} \ell(x) + \frac{1}{2} \ell(x') - \frac{1}{2c} \left(\ell(x) - \ell(x')\right)^2,
\]
where 
\[
c = \frac{\eta B^2}{4 \log\left(1 + \frac{\eta^2B^2}{2}\right)}.
\]

\subsection{Second claim}
Let $c>0$, we denote $\mathcal{E}(c)$ the set of functions $f: \mathcal{X} \to \mathbb{R}$, such that for any $x,x' \in \mathcal{X}$:

\begin{equation}\label{def:Ec}
f\left(\frac{x+x'}{2}\right) \le \frac{1}{2} f(x) + \frac{1}{2} f(x') - \frac{1}{2c} \left(f(x) - f(x')\right)^2.
\end{equation}
%Let $\mathcal{E} = \cup_{c>0} \mathcal{E}(c)$. Lemma~\ref{lem:assumption1} shows that any bounded and exp-concave function is in $\mathcal{E}$. Conversely, for any function $f$ in $\mathcal{E}(c)$ for some $c>0$, lemma~\ref{lem:recip1} shows that $f$ is bounded and lemma~\ref{lem:recip2} shows that $f$ is exp-concave. We conclude that the functions set $\mathcal{E}$ introduced above, coincides with the set of bounded and exp-concave functions. Furthermore, lemma~\ref{lem:recip1} shows that $\mathcal{E}$ contains the set of Lipschitz and strongly convex functions. This latter inclusion is strict : take for example the function $f: [-1,1] \to \mathbb{R}$ defined by $f(x) = x^4$, while $f \in \mathcal{E }$, $f$ is not strongly convex.
\begin{lemma}\label{lem:recip1}
	For any $c>0$, we have for any $f \in \mathcal{E}(c)$ 
	\[
	\sup_{x,x' \in \cX} \abs{f(x) - f(x')} \le c.
	\]
\end{lemma}
\begin{proof}
	Put $\Delta_{xx'} = f(x')-f(x)$,
	and $\Delta^* = \sup_{x,x' \in \cX} \Delta_{xx'}$. We first prove that $\Delta^* \leq 3c$.
	Assume this is not the case and let $x,x'\in \cX$ be such that $\Delta_{xx'} > 3c$. Let $z:=\frac{1}{2}\paren{x+x'}$. Using $f \in \cE(c)$, we obtain
	\[
	\Delta_{xz} = f(z) - f(x) \leq \frac{1}{2}\paren{f(x')-f(x)} - \frac{1}{2c}(f(x')-f(x))^2 = \frac{1}{2}\Delta_{xx'} - \frac{1}{2c}
	\Delta_{xx'}^2 \leq - \Delta_{xx'},
	\]
	where the last inequality holds because $\Delta_{xx'} >3c$. Hence $\Delta_{zx} >3c$ and in turn,
	if $x_1:=\frac{1}{2}\paren{x+z}$, reiterating the above argument we get $\Delta_{x_1z}>3c$ and
	in particular $f(x_1)<f(z)$. Also, we have $\Delta_{zx'} = \Delta_{zx} + \Delta_{xx'} > 3c$,
	therefore putting $x'_1:=\frac{1}{2}\paren{x'+z}$, again by the same token we get $f(x'_1)<f(z)$.
	This is a contradiction, since $z=\frac{1}{2}\paren{x_1+x'_1}$, thus Assumption~\ref{assump}
	implies that $f(z)\leq \max(f(x_1),f(x'_1))$.
	
	Since $\Delta^*$ is finite, $m:=\inf_{x \in X} f(x)$ is finite. For any $\eps>0$, let $x_\eps$ be
	such that $f(x_\eps) \leq m + \eps$. For any $x' \in X$, putting again $z:=\frac{1}{2}\paren{x+x'}$,
	it must be the case that $\Delta_{x_\eps z} \geq -\eps$, and using again the above display it must hold
	$-\eps \leq \Delta_{x_\eps z} \leq \frac{1}{2}\Delta_{x_\eps x'} - \frac{1}{2c}
	\Delta_{x_\eps x'}^2$. This implies $\Delta_{x_\eps x'} \leq c + G(\eps)$ for any $x' \in \cX$, with
	$G(\eps)=O(\eps)$. Since $\Delta^* \leq \eps + \sup_{x' \in \cX}\Delta_{x_\eps x'} $,
	we conclude to $\Delta^*\leq c$ by letting $\eps \rightarrow 0$.
\end{proof}

\begin{lemma}\label{lem:recip2}
	For any $c>0$, we have for any continuous function $f \in \mathcal{E}(c)$: $f$ is $(4/c)$-exp-concave.
\end{lemma}
\begin{proof}
	Fix $c>0$ and $f \in \mathcal{E}(c)$. Let $x, x' \in \cX$. Let us prove that
	\begin{equation}\label{eq:goal}
	\frac{1}{2} \exp\left\lbrace -\frac{4}{c} f(x) \right\rbrace + \frac{1}{2} \exp\left\lbrace - \frac{4}{c} f(x') \right \rbrace \le \exp\left\lbrace -\frac{4}{c} f\left(\frac{x+x'}{2}\right) \right\rbrace.
	\end{equation}
	
	Recall that since $f \in \mathcal{E}(c)$, inequality \eqref{def:Ec} gives
	\begin{equation*}
	\frac{2}{c^2} \left(f(x) - f(x')\right)^2 \le \frac{2}{c} f(x) + \frac{2}{c} f(x') - \frac{4}{c} f\left(\frac{x+x'}{2}\right).
	\end{equation*}
	We introduce the $\exp$ function on both sides of the inequality and use the first result of Lemma~\ref{lem:colg} to lower bound the left hand side. We have
	\begin{equation*}
	\frac{1}{2} \exp\left\lbrace \frac{2}{c}\left(f(x) - f(x')\right)\right\rbrace + \frac{1}{2} \exp\left\lbrace \frac{2}{c}\left(f(x') - f(x)\right)\right\rbrace\le \exp\left\lbrace  \frac{2}{c} f(x) + \frac{2}{c} f(x') \right\rbrace~\exp\left\lbrace - \frac{4}{c} f\left(\frac{x+x'}{2}\right) \right\rbrace,
	\end{equation*}
	which proves \eqref{eq:goal}. We conclude using the characterization provided by Lemma~\ref{lem:conv}.
\end{proof}
\subsection{Third claim}
\begin{lemma}\label{lem:recip3}
	Let $f: \mathcal{X} \to \mathbb{R}$ be a $L$-Lipschitz and $\rho$-strongly convex function, then $f \in \mathcal{E}\left(4L^2/\rho\right)$.
\end{lemma}
\begin{proof}
	
	By strong convexity of $f$, we have for any $x,x' \in \mathcal{X}$
	\begin{equation*}
	f\left(\frac{x+x'}{2}\right) \le \frac{1}{2} f(x) + \frac{1}{2} f(x') - \frac{\rho}{8} \norm{x-x'}^2.
	\end{equation*}
	Moreover, $f(.)$ is $L$-Lipschitz, hence: $\lvert f(x) - f(x')\rvert \le L \norm{x-x'}$. Therefore
	\begin{equation*}
	f\left(\frac{x+x'}{2}\right) \le \frac{1}{2} f(x) + \frac{1}{2} f(x') - \frac{\rho}{8L^2} \left(f(x) - f(x')\right)^2.
	\end{equation*}
\end{proof}

\section{Concentration inequality for martingales}

We recall Bennett's inequality:
\begin{theorem}\label{thm:bennett}
	Let $Z,Z_1, \dots, Z_n$ be i.i.d random variables with values in $[-B,B]$ and let $\delta>0$. Then with probability at least $1-\delta$ in $(Z_1, \dots, Z_n)$ we have
	\[
	\abs{\mathbb{E}[Z] - \frac{1}{n} \sum_{i=1}^{n} Z_i} \le \sqrt{\frac{2 \Var[Z] \log(2/\delta)}{n}} + \frac{2B\log(2/\delta)}{3n}.
	\]
\end{theorem}

We recall Freedman's inequality (the exposition here is lifted from~\citealp{fan2015exponential}).
Let $(\xi_i,\cF_i)_{i\geq 1}$ be a (super)martingale difference sequence. Define
$S_n:=\sum_{i=1}^n \xi_i$ (then $(S_n,\cF_n)$ is a (super)martingale), and
$\inner{S}_n:= \sum_{i=1}^n \e{\xi_i^2|\cF_{i-1}}$ the quadratic characteristic of $S$.

\begin{theorem}[Freedman's inequality]
	Assume $\xi_i\leq B$ for all $i\geq 1$, where $B$ is a constant. Then for all $t,v>0$:
	\begin{equation}
	\label{eq:freedman}
	\prob{S_k \geq t \text{ and } \inner{S}_k \leq v^2 \text{ for some } k\geq 1} \leq
	\exp\paren{-\frac{t^2}{2(v^2+Bt)}}.
	\end{equation}
\end{theorem}
The following direct consequence also appears in \citep[Lemma 3]{kakade2008generalization} for fixed $k$. Here we give a version that holds uniformly in $k$. See also \citep[Theorem 12]{gaillard2014second} for a related result.
\begin{corollary}\label{cor:concentration_mart}
	Assume $\xi_i\leq B$ for all $i\geq 1$, where $B$ is a constant. Then for all $\delta \in (0,1/3)$, with probability
	at least $1-3\delta$ it holds 
	\[
	\forall k \geq 1 : S_k \leq 2  \sqrt{\inner{S}_k\eps(\delta,k)}
	+ 4B \eps(\delta,k), 
	\]
	where $\eps(\delta,k) := \log \delta^{-1} + 2 \log(1 + \log_2^+ (\inner{S}_k /B^2))$.
	
	If $\abs{\xi_i} \leq B$ for all $i\geq 1$, observe that $\eps(\delta,k) \leq \log \delta^{-1} + O(\log \log k)$.
\end{corollary}
\begin{proof}
	By standard calculations, it holds that if $t\geq v \sqrt{2\log \delta^{-1}} + 2B \log \delta^{-1}$,
	then $\frac{t^2}{2(v^2+Bt)} \geq \log \delta^{-1}$. Therefore~\eqref{eq:freedman} implies
	that for any $v>0$ and $\delta \in (0,1)$, it holds
	\begin{equation}
	\label{eq:freedman2}
	\prob{\exists k\geq 1: S_k \geq \sqrt{2v^2\log \delta^{-1}} + 2B \log \delta^{-1} \text{ and } \inner{S}_k \leq v^2} \leq \delta.
	\end{equation}
	Denote $v_j^2 := 2^{j}B^2$, $\delta_j := (j \vee 1)^{-2} \delta$, $j\geq 0$, and define the non-decreasing sequence of stopping times $\tau_{-1}=1$ and $\tau_j := \min \set[1]{k \geq 1: \inner{S}_k > v_j^2}$ for $j\geq 0$. 
	Define the events for $j\geq 0$:
	\begin{align*}
	A_j & := \set{ \exists k\geq 1: S_k \geq \sqrt{2v_j^2\log \delta_j^{-1}} + 2B \log \delta_j^{-1} \text{ and } \inner{S}_k \leq v_j^2},\\
	A'_j & := \set{ \exists k \text{ with }  \tau_{j-1} \leq k < \tau_{j} : S_k \geq 2\sqrt{\inner{S}_k \eps(\delta,k)} + 4B \eps(\delta,k)}.
	\end{align*}
	From the definition of $v_j^2,\delta_j$, we have $j = \log_2 (v_j^2/B^2)$ for $j\geq 1$.
	For $j\geq 1$, $\tau_{j-1} \leq k < \tau_{j}$ implies $v_{j-1}^2 = v^2_j/2 < \inner{S}_k \leq v_{j}^2 $, and further
	\[
	\log \delta_j^{-1} = \log \delta^{-1} + 2 \log \log_2 (v_j^2 /B^2)
	% \leq \log \delta^{-1} + 2 \log(1 + \log_2 (\inner{S}_k /B^2)_+) =
	\leq \eps(\delta,k).\]
	Therefore it holds $A'_j \subseteq A_j$.
	Furthermore, for $j=0$, we have $v_0^2=B^2,\delta_0=\delta$.
	Further, if $k <\tau_0$ it implies $\inner{S}_k <B^2$ and therefore
	$\eps(\delta,k) = \log \delta^{-1}$.
	Thus, provided $\log \delta^{-1}\geq 1$ i.e. $\delta \leq 1/e$, it holds %$k<\tau_j$ implies $\inner{s}_k \leq B^2$, hence
	\begin{multline*}
	A'_0 \subseteq \set{ \exists k \text{ with } k<\tau_0 : S_k \geq 4B \log \delta_0^{-1}}\\
	\subseteq
	\set{ \exists k \geq 1 : S_k \geq  \sqrt{2v_0^2\log \delta_0^{-1}} + 2B \log \delta_0^{-1}
		\text{ and } \inner{S}_k \leq v_0^2 } = A_0.
	\end{multline*}
	Therefore, since by~\eqref{eq:freedman2} it holds $\prob{A_j} \leq \delta_j$ for all $j \geq 0$:
	\begin{align*}
	\prob{\exists k\leq n: S_k \geq 2  \sqrt{\inner{S}_k\eps(\delta,k) } + 4B \eps(\delta,k)}  
	= \prob[4]{\bigcup_{j \geq 0 } A'_j}
	\leq \prob[4]{\bigcup_{j \geq 0 } A_j}
	\leq \delta \sum_{j\geq 0} (j \vee 1)^{-2} \leq 3 \delta. 
	\end{align*}
	
\end{proof}

\begin{corollary}\label{cor:concentration_mart2}
	Assume $\xi_i \le b$ for all $i \ge 1$, where $b$ is a constant. Let $(\nu_t)_t$ denote an $\mathcal{F}_t$-measurable sequence, such that for any $k\ge 1$: $\langle S \rangle_k \le \sum_{i=1}^{k} \nu_i$. Then for all $c>0$ and $\delta \in (0, 1/3)$, with probability at least $1-3\delta$ it holds  
	\begin{equation*}
	\forall k \geq 1 : S_k - \frac{c}{b} \sum_{i=1}^{k} \nu_k  \leq  \left(\frac{8}{c}+4\right) \left( \log(\delta^{-1})+2\log_2^+\left(\frac{32+16c}{c^2} \right)\right) b.
	\end{equation*}
\end{corollary}
\begin{proof}
	Let $c>0$ and fix $\delta \in (0,1/3)$, we have using Corollary~\ref{cor:concentration_mart}: with probability at least $1-3\delta$, it holds for any $k \ge 1$
	\begin{align*}
	S_k - \frac{c}{b} \sum_{i=1}^{k} \nu_i &\le 2\sqrt{\langle S \rangle_k \epsilon(\delta,k)} + 4b \epsilon(\delta, k) - \frac{c}{b} \sum_{i=1}^{k} \nu_i \\
	&\le 2\sqrt{\langle S \rangle_k \epsilon(\delta,k)} + 4b \epsilon(\delta, k) - \frac{c}{b} \langle S \rangle_k\\
	&\le 2 \left(\frac{c}{4b}\langle S \rangle_k + \frac{4b}{c}\epsilon(\delta, k) \right) +4b \epsilon(\delta, k) - \frac{c}{b} \langle S \rangle_k\\
	&\le \left(\frac{8}{c} +4\right) b\epsilon(\delta,k) - \frac{c}{2b} \langle S \rangle_k\\
	&= \left(\frac{8}{c} +4\right) b \left(\log \delta^{-1} + 2 \log\left(1+ \log_2^+(\langle S \rangle_k / b^2)\right)\right) - \frac{c}{2b} \langle S \rangle_k\\
	&\le \left(\frac{8}{c} +4\right) b \left(\log \delta^{-1} + 2 \log_2^+(\langle S \rangle_k / b^2)\right) - \frac{c}{2b} \langle S \rangle_k
	\end{align*}
	The result follows by upper-bounding the function $x \to \log_2^+(x) - x/y$, for $x,y>0$ using Lemma~\ref{lem:pure_cal}.
\end{proof}

\section{Additional technical results}\label{sec:add_res}

The following lemma is a consequence of Corollary~\ref{cor:bennett}, the chaining rule (i.e cancellation in the sum of logarithmic terms) and Fubini's theorem. Let $(\hat{h}_{i,t})_{t \in \intr{T}, i \in \intr{K}}$ be a $\mathcal{F}_t$-adapted process.

For each $i\in \intr{K}$ and $t \in \intr{T}$ we define: $\hat{H}_{i,t} := \sum_{i=1}^{t} \hat{h}_{i,s}$, we use the convention that $\hat{H}_{i,0} = 0$. Let $t \in \intr{T}$ and $\lambda >0$, we define the sequence $(\hat{p}_{i,t})_{i \in \intr{K}}$:
\begin{equation}\label{eq:gibs_weights}
\hat{p}_{i,t} := \frac{\exp\left\{-\lambda \hat{H}_{i,t-1}\right\}}{\sum_{j=1}^{K}\exp\left\{-\lambda \hat{H}_{j,t-1}\right\}}.
\end{equation}
%For each $t \in \intr{T}$, let $\hat{Z}_{t} := \sum_{i=1}^{K} \exp\{-\lambda \hat{H}_{i,t}\}$ and define $M_t := \log\left(\hat{Z}_t\right) - \mathbb{E}_{t-1} \left[\log(\hat{Z}_t)\right]$.

For each $t \in \intr{T}$, define:
\begin{align}
\hat{Z}_{t} &:= \sum_{i=1}^{K} \exp\{-\lambda \hat{H}_{i,t}\} \label{eq:Z_t}\\
M_t &:= \log\left(\hat{Z}_t\right) - \mathbb{E}_{t-1} \left[\log(\hat{Z}_t)\right]. \label{eq:M_t}
\end{align}

\begin{lemma}\label{lem:log_tel}
	Let $b>0$ and $(\hat{h}_{i,t})_{t \in \intr{T}, i\in \intr{K}}$ be a sequence of numbers taking values in an interval of length $b$. For each $i\in \intr{K}$ and $t\in \intr{T}$, let $\mathbb{E}_{t-1}[\hat{h}_{i,t}] = h_{i,t}$.
	Let $(\alpha_t)_{t \in \intr{T}}$ be a sequence such that $\alpha_t$ is $\mathcal{F}_{t-1}$-measurable and: 
	\[
	\forall i \in \intr{K}, t\in \intr{T},~\abs{\hat{h}_{i,t} - \alpha_t} \le b. 
	\] 
	\noindent Then for any $\lambda \in \left( 0, 1/b\right)$, for all $t \in \intr{T}$ we have:
	\begin{equation*}
	\sum_{t=1}^{T}\sum_{i=1}^{K} \hat{p}_{i,t}~h_{i,t} \le \min_{i\in \intr{K}} \sum_{t=1}^{T}\hat{h}_{i,t} + \frac{\log(K)}{\lambda} + \frac{1}{\lambda} \sum_{t=1}^{T-1} M_t +  \lambda \sum_{t=1}^{T}\sum_{i=1}^{K} \hat{p}_{i,t}~\mathbb{E}_{t-1}\left[ \left(\hat{h}_{i,t} - \alpha_t\right)^2 \right],	
	\end{equation*}
	where the sequence $(\hat{p}_{i,t})_{t \in \intr{T}, i \in \intr{K}}$ is defined by \eqref{eq:gibs_weights} and $(M_t)$ is defined by \eqref{eq:M_t}.
\end{lemma}

\begin{proof}
	Let  $t\in \intr{T}$, we denote by $\hat{p}_t$ the probability distribution on $\intr{K}$ defined by the weights $(\hat{p}_{i,t})_{i \in \intr{K}}$. We apply Corollary~\ref{cor:bennett} to the random variable $X_t := \hat{h}_{I,t} - \alpha_t$, where $I$ is drawn from $\intr{K}$ following $\hat{p}_t$: for any $\lambda \in \left(0, 1/b\right)$,
	\begin{equation*}
	\log\left( \sum_{i=1}^{K} \hat{p}_{i,t} \exp\left\lbrace - \lambda \left(\hat{h}_{i,t} - \alpha_t\right)  \right \rbrace\right) \le -\lambda \sum_{i=1}^{K} \hat{p}_{i,t}\left(\hat{h}_{i,t} - \alpha_t\right) + \lambda^2 \sum_{i=1}^{K} \hat{p}_{i,t} \left(\hat{h}_{i,t} - \alpha_t\right)^2.
	\end{equation*}
	Rearranging terms we obtain:
	\begin{align*}
	\sum_{i=1}^{K} \hat{p}_{i,t}~\hat{h}_{i,t} &\le  \alpha_t - \frac{1}{\lambda} \log\left( \left(\sum_{i=1}^{K} \hat{p}_{i,t} \exp\{ -\lambda \hat{h}_{i,t} \} \right)  \exp\{ \lambda \alpha_t \}  \right) + \lambda \sum_{i=1}^{K} \hat{p}_{i,t} \left(\hat{h}_{i,t} - \alpha_t\right)^2 \\
	&= -\frac{1}{\lambda} \log\left( \sum_{i=1}^{K} \hat{p}_{i,t} \exp\{ -\lambda \hat{h}_{i,t} \} \right) + \lambda \sum_{i=1}^{K} \hat{p}_{i,t} \left(\hat{h}_{i,t} - \alpha_t\right)^2 \\
	&= -\frac{1}{\lambda} \left( \log\left(\hat{Z}_{t}\right) - \log\left(\hat{Z}_{t-1}\right)\right) + \lambda \sum_{i=1}^{K} \hat{p}_{i,t} \left(\hat{h}_{i,t} - \alpha_t\right)^2,
	\end{align*} 
	where $\hat{Z}_t$ is defined by \eqref{eq:Z_t}.
	Taking the conditional expectation with respect to $\mathcal{F}_{t-1}$ gives
	\begin{equation*}
	\sum_{i=1}^{K} \hat{p}_{i,t} h_{i,t} \le -\frac{1}{\lambda} \left( \mathbb{E}_{t-1}\left[ \log\left(\hat{Z}_{t}\right) \right] - \log\left(\hat{Z}_{t-1}\right)\right) + \lambda \sum_{i=1}^{K} \hat{p}_{i,t} \mathbb{E}_{t-1}\left[ \left(\hat{h}_{i,t} - \alpha_t\right)^2 \right].
	\end{equation*}
	
	\noindent Summing over $t \in \intr{T}$ we obtain:
	\begin{equation*}
	\sum_{t=1}^{T}\sum_{i=1}^{K} \hat{p}_{i,t}~h_{i,t} \le \frac{\log\left(Z_0\right)}{\lambda} - \frac{\log\left(\hat{Z}_T\right)}{\lambda} + \frac{1}{\lambda} \sum_{t=1}^{T-1} M_t +  \lambda \sum_{t=1}^{T}\sum_{i=1}^{K} \hat{p}_{i,t}~\mathbb{E}_{t-1}\left[ \left(\hat{h}_{i,t} - \alpha_t\right)^2 \right].
	\end{equation*}

	\noindent Finally observe that $Z_0 = K$ and that:
	\begin{align*}
	-\frac{1}{\lambda} \log\left(\hat{Z}_T\right) &=  -\frac{1}{\lambda} \log\left( \sum_{i} \exp\{-\lambda \hat{H}_{i,t} \} \right) \\
	&\le \min_{i \in \intr{K}} \hat{H}_{i,t}.
	\end{align*}
\end{proof}

\section{A preliminary result for the proof of Theorem~\ref{thm:1} and~\ref{thm:1bis}}
In this section we present two key results for the proof of Theorem~\ref{thm:1} and~\ref{thm:1bis}. Lemma~\ref{lem:keyf1} provides a bound for the cases $(p=2, m\ge 3)$ and $(p=2, m=2, \text{IC}=\text{False})$. Lemma~\ref{lem:keyf2} presents a similar bound for the particular case $(p=2, m = 2, \text{IC}=\text{True})$. We decided to separate these two settings because each one requires a different condition on $\lambda$.

We consider the notation of Algorithms~\ref{algo:1} and~\ref{algo:2}. 
In Algorithm~\ref{algo:1} ($m\ge3$), we take $A_t = I_t$. Recall that  $\tilde{m} = \max\{1, m-2\}$ (as defined in Section~\ref{sec:not}).
\begin{lemma}\label{lem:powerk}
	For any $k \ge 1$, 
	\[
	\mathbb{E}_{t-1}\left[\left(\hat{\ell}_{i,t} - \ell_{A_t,t}\right)^k\right] = \left( \frac{K}{\tilde{m}} \right)^{k-1} \mathbb{E}_{t-1}\left[\left(\ell_{i,t} - \ell_{A_t,t}\right)^k \right],
	\]
	where $\tilde{m} = \max\{1, m-2\}$.
\end{lemma}	
\begin{proof}
	Suppose that $m\ge 3$. Consider the notation of Algorithm~\ref{algo:1}. Let $k\ge1$, we have
	\begin{align*}
	\mathbb{E}_{t-1}\left[\left(\hat{\ell}_{i,t} - \ell_{A_t,t}\right)^k\right] &= \mathbb{E}_{t-1}\left[ \left(\frac{K}{m-2} \mathds{1}\left(i \in \cU_t\right) \ell_{i,t} + \left(1 - \frac{K}{m-2} \mathds{1}\left(i \in \cU_t\right)\right) \ell_{A_t,t} - \ell_{A_t,t}\right)^k \right]\\
	&= \mathbb{E}_{t-1}\left[ \left(\frac{K}{m-2} \mathds{1}\left(i \in \cU_t\right) \ell_{i,t} - \frac{K}{m-2} \mathds{1}\left(i \in \cU_t\right) \ell_{A_t,t} \right)^k \right]\\
	&= \left( \frac{K}{m-2}\right)^k \mathbb{E}_{t-1}\left[\mathds{1}\left(i \in \cU_t\right)\right] \left(\ell_{i,t} - \ell_{A_t,t}\right)^k\\
	&= \left( \frac{K}{m-2} \right)^{k-1} \mathbb{E}_{t-1}\left[\left(\ell_{i,t} - \ell_{A_t,t}\right)^k \right],
	\end{align*}
	where we used the fact that $U_t$ and $A_t$ are independent conditionally to $\cF_{t-1}$. 
	
	\noindent Suppose that $m=2$. Consider the notation of Algorithm~\ref{algo:2}. Let $k \ge 1$, we have
	\begin{align*}
	\mathbb{E}_{t-1}\left[\left(\hat{\ell}_{i,t} - \ell_{A_t,t}\right)^k\right] &= \mathbb{E}_{t-1}\left[\left(\hat{\ell}_{i,t} - \ell_{A_t,t}\right)^k\right]\\
	&= \mathbb{E}_{t-1}\left[ \left(K \mathds{1}\left(B_t = i\right) \ell_{i,t} + \left(1 - K \mathds{1}\left(B_t = i\right)\right) \ell_{A_t,t} - \ell_{A_t,t}\right)^k \right]\\
	&= K^k \mathbb{E}_{t-1}\left[ \mathds{1}\left(B_t = i\right) \left(\ell_{i,t} - \ell_{A_t,t}\right)^k\right]\\
	&= K^{k-1} \mathbb{E}_{t-1}\left[ \left(\ell_{i,t} - \ell_{A_t,t}\right)^k \right].
	\end{align*}
\end{proof}

\noindent Introduce the notation
\begin{align}
\hat{\mu}_t &:= \sum_{i \in \intr{K}} \hat{p}_{i,t} \ell_{i,t},\label{eq:def_mu_1}\\
\hat{\xi}_t &:= \frac{1}{2}\sum_{i,j \in \intr{K}} \hat{p}_{i,t}\hat{p}_{j,t}~(\ell_{i,t} - \ell_{j,t})^2,\label{eq:def_xi_1}
\end{align}
where $(\hat{p}_{i,t})$ is defined in \eqref{eq:gibs_weights}.
\noindent For each $t\in \intr{T}$, let
\begin{align}
\hat{Z}_t &= \sum_{i=1}^{K} \exp\left\lbrace -\lambda \hat{L}_{i,t} + \lambda^2 \hat{V}_{i,t}\right\rbrace\notag\\
M_t &= \log\left(\hat{Z}_t\right) - \mathbb{E}_{t-1}\left[\hat{Z}_t\right] \label{eq:def_mt},
\end{align}
where $\hat{L}_{i,t} = \sum_{s=1}^{t} \hat{\ell}_{i,t}$ and $\hat{V}_{i,t} = \sum_{s=1}^{t} \hat{v}_{i,t}$, in agreement with the notation used in Algorithms~\ref{algo:1} and~\ref{algo:2}, and in Section~\ref{sec:add_res}.

\begin{lemma}\label{lem:key0}
	Let $\lambda \in \left(0, \frac{2\tilde{m}}{K}\bar{\lambda}\right)$, where $\bar{\lambda}$ is defined in \eqref{eq:c_cst} and $\tilde{m} = \max\{m-2, 1\}$. For each $i \in \intr{K}$, $t\in \intr{T}$, let $\hat{h}_{i,t} = \hat{\ell}_{i,t}-\lambda \hat{v}_{i,t}$. We have
	\[
	\sum_{t=1}^{T} \hat{\mu}_t \le \min_{i \in \intr{K}} \sum_{t=1}^{T} \hat{h}_{i,t} + \frac{1}{\lambda} \sum_{t=1}^{T} M_t + \frac{\log(K)}{\lambda} + \frac{11\lambda K}{\tilde{m}} \sum_{t=1}^{T} \hat{\xi}_t,
	\]
	where $\hat{\mu}_t$ is defined in \eqref{eq:def_mu_1}, $\hat{\xi}_t$ is defined in \eqref{eq:def_xi_1} and $M_t$ is defined in \eqref{eq:def_mt}.
\end{lemma}
\begin{proof}
	Let $h_{i,t} := \mathbb{E}_{t-1}[\hat{h}_{i,t}] = \ell_{i,t} - \lambda \mathbb{E}_{t-1}\left[\hat{v}_{i,t}\right] $, we apply Lemma~\ref{lem:log_tel} to the sequence $(\hat{h}_{i,t})_{i,t}$. We take $\alpha_t = \hat{\mu}_t$, which is an $\cF_{t-1}$-measurable process. For each $i \in \intr{K}$ and $t \ge 0$, we have
	\begin{equation}\label{eq:h} 
	\sum_{t=1}^{T} \sum_{i=1}^{K} \hat{p}_{i,t} h_{i,t}  \le \min_{i\in \intr{K}} \sum_{t=1}^{T} \hat{h}_{i,t} + \frac{\log(K)}{\lambda} + \frac{1}{\lambda} \sum_{t=1}^{T} M_t + \lambda \sum_{t=1}^{T} \sum_{i=1}^{K} \hat{p}_{i,t} \mathbb{E}_{t-1}\left[ \left(\hat{h}_{i,t} - \hat{\mu}_t\right)^2\right].
	\end{equation}
	
	\noindent Now, let us develop a lower bound on the left hand side of the inequality above. Recall that in Algorithm~\ref{algo:1}, we take $A_t = I_t$, then $A_t\sim \hat{p}_t$. In Algorithm~\ref{algo:2}, Lemma~\ref{lem:sample1} shows that $A_t\sim\hat{p}_t$. Fix $t \in \intr{T}$, we have:
	\begin{align}
	\sum_{i=1}^{K} \hat{p}_{i,t} h_{i,t} &= \sum_{i =1}^{K} \hat{p}_{i,t} \left(\ell_{i,t} - \lambda \mathbb{E}_{t-1}[\hat{v}_{i,t}]\right) \notag\\
	&= \sum_{i=1}^{K} \hat{p}_{i,t}\ell_{i,t} - \lambda \sum_{i=1}^{K} \hat{p}_{i,t} \mathbb{E}_{t-1}\left[\left(\hat{\ell}_{i,t} - \ell_{A_t,t}\right)^2\right]\notag\\
	&= \sum_{i=1}^{K} \hat{p}_{i,t}\ell_{i,t} - \lambda \frac{K}{\tilde{m}} \left( \sum_{i=1}^{K} \hat{p}_{i,t} \left(\ell_{i,t} - \hat{\mu}_t\right)^2\right) -\lambda \frac{K}{\tilde{m}} \mathbb{E}_{t-1}\left[ \left(\ell_{A_t,t} - \hat{\mu}_t\right)^2\right]  \notag\\
	&= \hat{\mu}_t - 2 \lambda\frac{K}{\tilde{m}} \hat{\xi}_t,\label{eq:term2_1}
	\end{align} 
	where we used in the second line the definition $\hat{v}_{i,t} = \left(\hat{\ell}_{i,t} - \ell_{A_t,t}\right)^2$, Lemma~\ref{lem:powerk} with $k=2$ in the third line, and the fact that $A_t$ is distributed following $\hat{p}$ in the third and fourth line.
	
	Next, we develop an upper bound on the last term of the right hand side of \eqref{eq:h}. We have
	\begin{equation}\label{eq:term2_2}
	\sum_{t=1}^{T} \sum_{i=1}^{K} \hat{p}_{i,t}~\mathbb{E}_{t-1}\left[ \left(\hat{h}_{i,t} - \hat{\mu}_t\right)^2\right] \le 2 \sum_{t=1}^{T} \sum_{i=1}^{K} \hat{p}_{i,t} \left\lbrace\mathbb{E}_{t-1}\left[ \left(\hat{\ell}_{i,t} - \hat{\mu}_t\right)^2\right] + \lambda^2 \mathbb{E}_{t-1}\left[ \hat{v}_{i,t}^2\right] \right\rbrace. 
	\end{equation}
	
	\noindent Fix $t \in \intr{T}$. Let us bound each of the terms in the right hand side of the inequality above
	\begin{align}
	\sum_{i=1}^{K} \hat{p}_{i,t} \mathbb{E}_{t-1}\left[ \left(\hat{\ell}_{i,t} - \hat{\mu}_t\right)^2\right] &\le  \sum_{i=1}^{K} 2\hat{p}_{i,t} \left( \mathbb{E}_{t-1}\left[ \left(\hat{\ell}_{i,t} - \ell_{A_t,t}\right)^2\right] + \mathbb{E}_{t-1}\left[ \left(\ell_{A_t,t} - \hat{\mu}_t\right)^2\right]\right)\notag\\ 
	&=  2 \mathbb{E}_{t-1}\left[ \left(\ell_{A_t,t} - \hat{\mu}_t\right)^2\right]+ 2\frac{K}{\tilde{m}} \sum_{i=1}^{K}\hat{p}_{i,t}\mathbb{E}_{t-1} \left[  \left( \ell_{i,t} - \ell_{A_t,t}\right)^2 \right] \notag\\
	&= 2 \hat{\xi}_t + 2\frac{K}{\tilde{m}} \sum_{i=1}^{K}\hat{p}_{i,t}\left\lbrace   \left( \ell_{i,t} - \hat{\mu}_t \right)^2 +  \mathbb{E}_{t-1} \left[  \left( \ell_{A_t,t} - \hat{\mu}_t \right)^2 \right]\right \rbrace\notag \\
	&\le  \frac{6K}{\tilde{m}}  \hat{\xi}_t\label{eq:term2_3},
	\end{align}
	where we used Lemma~\ref{lem:powerk} for the second line. Moreover, using the same Lemma~\ref{lem:powerk} with $k=4$, we have
	\begin{align}
	\sum_{i=1}^{K} \hat{p}_{i,t}  \mathbb{E}_{t-1}\left[ \hat{v}_{i,t}^2\right] &= \sum_{i=1}^{K} \hat{p}_{i,t} \left(\frac{K}{\tilde{m}}\right)^3 \mathbb{E}_{t-1}\left[ \left(\ell_{i,t} - \ell_{A_t,t}\right)^4\right] \notag\\
	&\le \left(\frac{K}{\tilde{m}}\right)^3 B^2 \sum_{i=1}^{K} \hat{p}_{i,t} \mathbb{E}_{t-1}\left[ \left(\ell_{i,t} - \ell_{A_t,t}\right)^2\right] \notag\\
	&= 2 \left(\frac{K}{\tilde{m}}\right)^3 B^2 \hat{\xi}_t\label{eq:term2_4}.
	\end{align}
	
	\noindent We plug the bounds obtained from \eqref{eq:term2_3} and \eqref{eq:term2_4} into inequality \eqref{eq:term2_1}, and obtain
	\begin{equation}\label{eq:inter1}
	\sum_{t=1}^{T} \sum_{i=1}^{K} \hat{p}_{i,t}~\mathbb{E}_{t-1}\left[ \left(\hat{h}_{i,t} - \hat{\mu}_t\right)^2\right] \le 2\left(\frac{6K}{\tilde{m}} + 2\lambda^2\frac{K^3}{(\tilde{m})^3}B^2\right) \sum_{t=1}^{T} \hat{\xi_t}.
	\end{equation}
	Recall that by definition \eqref{eq:c_cst}, $\bar{\lambda} \le \frac{1}{B}$. Hence, $\lambda<\frac{2\tilde{m}}{K}\bar{\lambda}$ gives 
	\[
	\lambda^2\frac{K^2}{\tilde{m}^2}B^2 \le 4,
	\]
	\noindent we plug this bound into \eqref{eq:inter1} and obtain 
	\begin{equation}\label{eq:term2_5}
	\sum_{t=1}^{T} \sum_{i=1}^{K} \hat{p}_{i,t}~\mathbb{E}_{t-1}\left[ \left(\hat{h}_{i,t} - \hat{\mu}_t\right)^2\right] \le 20\frac{K}{\tilde{m}} \sum_{t=1}^{T} \hat{\xi_t}.
	\end{equation}
	\noindent Next, we plug the bounds obtained in \eqref{eq:term2_1} and \eqref{eq:term2_5} into \eqref{eq:h} to obtain
	\begin{equation*}
	\sum_{t=1}^{T}\hat{\mu}_t \le \min_{i\in \intr{K}} \sum_{t=1}^{T} \hat{h}_{i,t} + \frac{1}{\lambda} \sum_{t=1}^{T}M_t + \frac{\log(K)}{\lambda} + \frac{22\lambda K}{\tilde{m}} \sum_{t=1}^{T} \hat{\xi}_t.
	\end{equation*}
	
\end{proof}

\begin{lemma}\label{lem:mt}
	Let $\lambda \in \left(0, \frac{2\tilde{m}}{K}\bar{\lambda}\right)$, where $\bar{\lambda}$ is defined in \eqref{eq:c_cst} and $\tilde{m} = \max\{1, m-2\}$. Consider the martingale difference sequence $(M_t)_{t\in \intr{T}}$ defined in \eqref{eq:def_mt}. We have
	\begin{itemize}
		\item $\forall t\in \intr{T}: \abs{M_t} \le 3\lambda \frac{K}{\tilde{m}} B$.
		\item $\sum_{t=1}^{T} \mathbb{E}\left[M_t^2\right] \le 5\frac{K}{\tilde{m}}\lambda^2 \sum_{t=1}^{T}\hat{\xi}_t$.
	\end{itemize}		
\end{lemma}
\begin{proof}
	Observe that the sequence $\left(M_t, \mathcal{F}_t\right)_{t \in \intr{T}}$ is a martingale difference.
	For any $t \in \intr{T}$, we have
	\begin{align*}
	M_t &= \mathbb{E}\left[\log\left(\hat{Z}_{t+1}\right) | \mathcal{F}_{t} \right] - \log\left(\hat{Z}_t\right)\\
	&= \log\left(\frac{\hat{Z}_t}{\hat{Z}_{t-1}}\right) - \mathbb{E}_{t-1}\left[\log\left(\frac{\hat{Z}_t}{\hat{Z}_{t-1}}\right)\right]\\
	&= \log\left( \sum_{i=1}^{K} \hat{p}_{i,t} \exp\{ -\lambda \hat{\ell}_{i,t} + \lambda^2 \hat{v}_{i,t} \}\right) - \mathbb{E}_{t-1}\left[ \log\left( \sum_{i=1}^{K} \hat{p}_{i,t} \exp\{ -\lambda \hat{\ell}_{i,t} + \lambda^2 \hat{v}_{i,t} \}\right) \right],
	\end{align*} 
	where we used the fact that $\hat{Z}_{t-1}$ is $\cF_{t-1}$-measurable in the second line.
	
	The loss function $\ell(.,y)$ is $B$-range-bounded for any $y$. Let $c_{\min}$ and $c_{max}$ denote the lower and upper bounds, respectively, for the values of $\ell$ ($c_{\max} - c_{\min} \le B$). 
	\noindent Therefore, for any $i \in \intr{K}$, $\hat{\ell}_{i,t} \in \left[c_{\min}-\frac{K}{\tilde{m}} B, c_{\max} + \frac{K}{\tilde{m}} B \right]$ and $\hat{v}_{i,t} \in [0, (\frac{K}{\tilde{m}})^2B^2]$. Therefore
	\begin{equation*}
	\exp\left(\lambda c_{\max}-\frac{K}{\tilde{m}} \lambda B \right) \le \exp\left(-\lambda \hat{\ell}_{i,t} + \lambda^2 \hat{v}_{i,t}\right) \le \exp\left(-\lambda c_{\min}+\lambda \frac{K}{\tilde{m}}B+2\lambda^2 \frac{K^2}{\tilde{m}^2}B^2\right).
	\end{equation*}
	Hence
	\begin{equation*}
	\lambda c_{\max} -\lambda \frac{KB}{\tilde{m}} \le \log\left( \sum_{i=1}^{K} \hat{p}_{i,t} \exp\{ -\lambda \hat{\ell}_{i,t} + \lambda^2 \hat{v}_{i,t} \} \right) \le -\lambda c_{\min} +\lambda \frac{KB}{\tilde{m}} + 2 \lambda^2 \frac{K^2B^2}{\tilde{m}^2}
	\end{equation*}
	\noindent Recall that $M_t$ is a centered variable and $\lambda < \frac{\tilde{m}}{128KB}$. Therefore
	
	\begin{equation}\label{eq:bound_mt}
	\abs{M_t} \le 4\lambda \frac{K}{\tilde{m}}B. 
	\end{equation}

	\noindent Now, let us bound the quadratic characteristic of $\left(M_t\right)_t$. We have
	\begin{align}
	\mathbb{E}_{t-1}\left[ M_t^2\right] &= \Var_{t-1}\left( \log\left(\hat{Z}_t\right)\right) \nonumber \\
	&= \Var_{t-1}\left( \log\left(\hat{Z}_t\right) - \log\left(\hat{Z}_{t-1}\right)\right), \label{eq:bound_Mt}
	\end{align}
	where we used the fact that $\hat{Z}_{t-1}$ is $\mathcal{F}_{t-1}$-measurable. 
	
	\noindent Furthermore we have
	\begin{align*}
	\hat{Z}_{t} &= \sum_{i=1}^K \exp\left(- \lambda \hat{L}_{i,t} + \lambda^2 \hat{V}_{i,t} \right)\\
	&= \sum_{i=1}^{K} \exp\left( - \lambda \hat{L}_{i,t-1} + \lambda^2 \hat{V}_{i,t}\right) \exp\left(- \lambda \hat{\ell}_{i,t} + \lambda^2 \hat{v}_{i,t}\right)\\
	&= \sum_{i=1}^{K} \hat{p}_{i,t} \hat{Z}_{t-1} \exp\left(- \lambda \hat{\ell}_{i,t} + \lambda^2 \hat{v}_{i,t}\right).
	\end{align*}
	Hence
	\begin{align}
	\frac{\hat{Z}_{t}}{\hat{Z}_{t-1}} &= \sum_{i=1}^{K} \hat{p}_{i,t} \exp\left(- \lambda \hat{\ell}_{i,t} + \lambda^2 \hat{v}_{i,t}\right) \notag\\
	&= \sum_{i=1}^{K} \hat{p}_{i,t} \exp\left(- \lambda \left( \ell_{A_t,t} + \frac{K}{\tilde{m}} \mathds{1}\left(i \in \cU_t\right) (\ell_{i,t} - \ell_{A_t,t})\right) + \lambda^2 \frac{K^2}{\tilde{m}^2} \mathds{1}\left(i \in \cU_t\right) (\ell_{i,t} - \ell_{A_t,t})^2  \right)\notag\\
	&= \exp\left(-\lambda \ell_{A_t,t} \right)  \sum_{i=1}^{K} \hat{p}_{i,t} \exp\left( -\lambda \frac{K}{\tilde{m}} \mathds{1}\left(i \in \cU_t\right) (\ell_{i,t} - \ell_{A_t,t}) + \lambda^2 \frac{K^2}{\tilde{m}^2} \mathds{1}\left(i \in \cU_t\right) (\ell_{i,t} - \ell_{A_t,t})^2  \right) \notag\\
	&= \exp\left(-\lambda \ell_{A_t,t} \right)  \mathbb{E}_{A'_t} \left[ \exp\left( -\lambda \frac{K}{\tilde{m}} \mathds{1}\left(A'_t \in \cU_t\right) (\ell_{A'_t,t} - \ell_{A_t,t}) + \lambda^2 \frac{K^2}{\tilde{m}^2} \mathds{1}\left(A'_t \in \cU_t\right) (\ell_{A'_t,t} - \ell_{A_t,t})^2 \right) \right] \label{eq:Mt_1},
	\end{align}
	
	\noindent where $A'_t$ is a random variable, independent of $A_t$, such that for each $i\in \intr{K}$, $\mathbb{P}\left(A'_t = i\right) = \hat{p}_{i,t}$, and $\mathbb{E}_{A'_t}$ is the expectation with respect to the random variable $A'_t$. So as to ease notation, denote
	\begin{equation*}
	D_t := \frac{K}{\tilde{m}} \mathds{1}\left(A'_t \in \cU_t\right) (\ell_{A'_t,t} - \ell_{A_t,t})- \lambda \frac{K^2}{\tilde{m}^2} \mathds{1}\left(A'_t \in \cU_t\right) (\ell_{A'_t,t} - \ell_{A_t,t})^2.
	\end{equation*}
	We take the logarithm of both sides of inequality \eqref{eq:Mt_1}, we have
	\begin{equation*}
	\log\left(\hat{Z}_{t}\right) - \log\left(\hat{Z}_{t-1}\right) = - \lambda \ell_{A_t,t} + \log\left( \mathbb{E}_{A'_t} \left[ \exp\left( -\lambda D_t  \right) \right] \right). 
	\end{equation*}
	
	\noindent We inject the equality above in \eqref{eq:bound_Mt}. We obtain
	\begin{align}
	\mathbb{E}_{t-1}\left[ M_t^2\right] &= \Var_{t-1}\left(- \lambda \ell_{A_t,t} + \log\left( \mathbb{E}_{A'_t} \left[ \exp\left( -\lambda D_t  \right) \right] \right)  \right)\nonumber\\
	&\le 2\Var_{t-1}\left( \lambda \ell_{A_t,t}\right) + 2 \Var_{t-1} \left( \log\left(\mathbb{E}_{A'_t}\left[\exp\left(-\lambda D_t\right)\right] \right)\right)\nonumber\\
	&\le 2\Var_{t-1}\left( \lambda \ell_{A_t,t}\right) + 2 \mathbb{E}_{t-1} \left[ \log^2\left(\mathbb{E}_{A'_t}\left[\exp\left(-\lambda D_t\right)\right] \right)\right]. \label{eq:quad_Mt1}
	\end{align}
	
	\noindent Observe that
	\begin{equation*}
	\abs{\lambda D_t} = \abs{\lambda \frac{K}{\tilde{m}} \mathds{1}\left(A'_t \in \cU_t\right)\left(\ell_{A'_t,t} - \ell_{A_t,t}\right)- \lambda^2 \frac{K^2}{\tilde{m}^2} \mathds{1}\left(A'_t \in \cU_t\right) (\ell_{A'_t,t} - \ell_{A_t,t})^2 } \le \frac{1}{5}.
	\end{equation*} 
	where we used $\lambda \in \left(0, \frac{\tilde{m}}{128KB}\right)$. 
	
	\noindent The function $x\mapsto \log^2(x)$ is convex on $[e^{-1},e]$. Hence, using Jensen's inequality, we have
	\begin{align}
	\mathbb{E}_{t-1} \left[ \log^2\left(\mathbb{E}_{A'_t}\left[\exp\left(-\lambda D_t\right)\right]\right)\right] &\le \mathbb{E}_{t-1}\mathbb{E}_{A'_t}\left[ \log^2\left(\exp\left(-\lambda D_t\right) \right)\right]\nonumber\\
	&= \mathbb{E}_{t-1}\mathbb{E}_{A'_t}\left[ \lambda^2 D_t^2\right]\label{eq:quad_Mt2}
	\end{align}
	
	\noindent From \eqref{eq:quad_Mt1} and \eqref{eq:quad_Mt2}, we conclude that
	\begin{align}
	\mathbb{E}_{t-1}\left[ M_t^2\right] &\le 2\lambda^2 \Var_{t-1}(\ell_{A_t,t}) + 2 \mathbb{E}_{t-1}\mathbb{E}_{A'_t}\left[ \lambda^2 D_t^2\right]\notag\\
	&\le 2\lambda^2 \hat{\xi}_t+2 \mathbb{E}_{t-1}\mathbb{E}_{A'_t}\left[ \lambda^2 D_t^2\right]\label{eq:inter2}.
	\end{align}
	where we used $\Var_{t-1}(\ell_{A_t,t}) = \hat{\xi}_t$. Furthermore:
	\begin{align}
	\mathbb{E}_{t-1}\mathbb{E}_{A'_t}\left[ \lambda^2 D_t^2\right] &\le 2 \mathbb{E}_{t-1}\mathbb{E}_{A'_t}\left[   \frac{\lambda^2K^2}{\tilde{m}^2}  \mathds{1}\left(A'_t \in \cU_t\right) \left(\ell_{A'_t,t} - \ell_{A_t,t}\right)^2 +  \frac{K^4\lambda^4}{\tilde{m}^4}  \mathds{1}\left(A'_t \in \cU_t\right) \left(\ell_{A'_t,t} - \ell_{A_t,t}\right)^4 \right]\notag\\
	&\le 2\left(\frac{\lambda^2 K^2}{\tilde{m}^2} + \frac{\lambda^4K^4}{\tilde{m}^4}B^2\right)  \mathbb{E}_{t-1}\mathbb{E}_{A'_t}\left[\mathds{1}\left(A'_t \in \cU_t\right) (\ell_{A'_t,t} - \ell_{A_t,t})^2 \right]\notag\\
	&\le 3\frac{\lambda^2 K^2}{\tilde{m}^2}   \mathbb{E}_{t-1}\mathbb{E}_{A'_t}\left[\mathds{1}\left(A'_t \in \cU_t\right) (\ell_{A'_t,t} - \ell_{A_t,t})^2 \right]\notag\\
	&\le 3\frac{\lambda^2 K^2}{\tilde{m}^2}  \mathbb{E}_{A'_t}\left[\mathbb{E}_{t-1}\left[\mathds{1}\left(A'_t \in \cU_t\right)\right] \mathbb{E}_{t-1}\left[(\ell_{A'_t,t} - \ell_{A_t,t})^2 \right]\right]\notag\\
	&=  3\frac{\lambda^2 K^2}{\tilde{m}^2} \frac{\tilde{m}}{K} \sum_{i,j=1}^{K} \hat{p}_{i,t}\hat{p}_{j,t} \left(\ell_{i,t} - \ell_{j,t}\right)^2\notag\\
	&=3 \frac{K}{\tilde{m}} \lambda^2 \hat{\xi}_t\label{eq:inter3},
	\end{align}
	where we used the independence of $U_t$ and $A_t$ conditionally to $\cF_{t-1}$.
	
	\noindent We plug \eqref{eq:inter3} into \eqref{eq:inter2}.
	Therefore, it holds
	\begin{align*}
	\sum_{t=1}^{T}\mathbb{E}_{t-1}\left[ M_t^2\right] &\le\sum_{t=1}^{T} \left( 2\lambda^2 \hat{\xi}_t + 3\frac{K}{\tilde{m}} \lambda^2 \hat{\xi}_t \right) \\
	&\le 5\frac{K}{\tilde{m}} \lambda^2 \sum_{t=1}^{T}\hat{\xi}_t.
	\end{align*}
\end{proof}

The following lemma provides a bound with high probability on the quantity $\hat{L}_{i,T} - \lambda \hat{V}_{i,T}$, for each $i \in \intr{K}$.

\begin{lemma}\label{lem:minL}
	For any $i \in \intr{K}$ and $\lambda \in (0, \frac{\tilde{m}\bar{\lambda}}{128 K })$, with $\bar{\lambda}$ defined in \eqref{eq:c_cst} and $\tilde{m} = \max\{1, m-2\}$. We have for any $\delta \in (0,1/3)$, with probability at least $1 - 6\delta$:
	\[
	\hat{L}_{i,T} - \lambda \hat{V}_{i,T} \le L_{i,T} +\frac{721}{\lambda} \log\left(\frac{\tilde{m}}{KB\lambda \delta}\right). 
	\] 
\end{lemma}
\begin{proof}
	Let $i \in \intr{K}$. Recall that we have for any $t \in \intr{T}$
	\begin{align*}
	\hat{\ell}_{i,t} - \ell_{i,t} &= \left(\frac{K}{\tilde{m}} \mathds{1}\left(i \in \cU_t\right) -1 \right) \left( \ell_{i,t} - \ell_{A_t,t}\right)\\
	\hat{\ell}_{i,t} - \ell_{A_t,t} &= \frac{K}{\tilde{m}} \mathds{1}\left(i \in \cU_t\right) \left(\ell_{i,t} - \ell_{A_t,t}\right).
	\end{align*}

	We introduce the following notation
	\begin{equation*}
	\nu_{i,t} := \mathbb{E}_{t-1}\left[\left(\ell_{i,t} - \ell_{A_t,t}\right)^2\right].
	\end{equation*}

	\noindent We have
	\begin{align}
	\hat{L}_{i,T} - \lambda \hat{V}_{i,T} &=L_{i,T} + \sum_{t=1}^{T} \left(\hat{\ell}_{i,t} - \ell_{i,t} \right)  - \lambda \sum_{t=1}^{T} \left(\frac{K}{\tilde{m}} \right)^2 \mathds{1}\left(i \in \cU_t\right) (\ell_{i,t} - \ell_{A_t,t})^2 \notag\\
	&= L_{i,T} + \underbrace{\sum_{t=1}^{T} \left(\hat{\ell}_{i,t} - \ell_{i,t}\right) - \lambda\frac{K}{2\tilde{m}}\sum_{t=1}^{T} \nu_{i,t}}_{\text{Term 21}}\notag\\
	&\qquad \qquad + \underbrace{\lambda\frac{K}{2\tilde{m}}\sum_{t=1}^{T} \nu_{i,t} -\lambda \sum_{t=1}^{T} \left(\frac{K}{\tilde{m}} \right)^2 \mathds{1}\left(i \in \cU_t\right) (\ell_{i,t} - \ell_{A_t,t})^2}_{\text{Term 22}}.\label{eq:term0}
	\end{align}
	
	\paragraph{Bounding Term 21:}
	Observe that $(\hat{\ell}_{i,t} - \ell_{i,t})_t$ is a martingale difference with respect to the filtration $\mathcal{F}$, bounded in absolute value by $\frac{K}{\tilde{m}}B$. 
	Let us bound its quadratic characteristic. Recall that $A_t$ and $\cU_t$ are independent conditionally to $\cF_{t-1}$. We have
	\begin{align*}
	\sum_{t=1}^{T} \mathbb{E}_{t-1} \left[ (\hat{\ell}_{i,t} - \ell_{i,t})^2 \right] &= \sum_{t=1}^{T} \mathbb{E}_{t-1}\left[ \left(1 - \frac{K}{\tilde{m}}\mathds{1}\left(i \in \cU_t\right)\right)^2 (\ell_{i,t} - \ell_{A_t,t})^2 \right]\\
	&= \sum_{t=1}^{T} \mathbb{E}_{t-1}\left[ \left(1 - \frac{K}{\tilde{m}}\mathds{1}\left(i \in \cU_t\right)\right)^2\right] \mathbb{E}_{t-1}\left[ (\ell_{i,t} - \ell_{A_t,t})^2 \right]\\
	&\le \frac{K}{\tilde{m}} \sum_{t=1}^{T} \mathbb{E}_{t-1} \left[(\ell_{i,t} - \ell_{A_t,t})^2\right]\\
	&= \frac{K}{\tilde{m}} \sum_{t=1}^{T} \nu_{i,t}.
	\end{align*}

	\noindent Next, we apply Corollary~\ref{cor:concentration_mart2} to the sequence $(\hat{\ell}_{i,t} - \ell_{i,t})_{t \in \intr{T}}$: We take $c=\lambda K B/(4\tilde{m}) \le 1$, with probability at least $1 - 3\delta$, it holds
	
	\begin{equation}\label{eq:term11}
	\sum_{t=1}^{T} \left( \hat{\ell}_{i,t} - \ell_{i,t} \right) - \lambda \frac{K}{2\tilde{m}} \sum_{t=1}^{T} \nu_{i,t}   \le  \frac{720}{\lambda} \log\left(\frac{\tilde{m}}{KB\lambda \delta}\right).
	\end{equation}
	
	\paragraph{Bounding Term 22:}
	Define the sequence $(Q_t)_{t\in \intr{T}}$ as follows:
	\begin{equation*}
	Q_t :=   -\lambda\frac{K^2}{\tilde{m}^2} \mathds{1}(i \in \cU_t) (\ell_{i,t} - \ell_{A_t,t})^2 + \lambda \frac{K}{\tilde{m}} \nu_{i,t} .
	\end{equation*}
	Notice that $(Q_t)$ is a martingale difference sequence with respect to the filtration $\mathcal{F}$, and bounded in absolute value by $2\lambda\frac{K^2B^2}{\tilde{m}^2}$. Let us bound its quadratic characteristic. We have
	
	\begin{align*}
	\sum_{t=1}^{T} \mathbb{E}_{t-1} \left[ Q_t^2 \right]&\le \lambda^2\sum_{t=1}^{T} \mathbb{E}_{t-1}\left[\frac{K^4}{\tilde{m}^4} \mathds{1}\left(i \in \cU_t\right) \left( \ell_{i,t} - \ell_{A_t,t} \right)^4 \right] \\
	&\le \lambda^2\frac{K^4B^2}{\tilde{m}^4}\sum_{t=1}^{T} \mathbb{E}_{t-1}\left[ \mathds{1}\left(i \in \cU_t\right)\right]\mathbb{E}_{t-1}\left[ \left( \ell_{i,t} - \ell_{A_t,t} \right)^2 \right]\\
	&= \frac{K^3\lambda^2B^2}{\tilde{m}^3}\sum_{t=1}^{T} \nu_{i,t}.
	\end{align*}
	
	\noindent Next, we apply Corollary~\ref{cor:concentration_mart2} to this sequence. We take $c=1$, we have with probability at least $1-3\delta$:
	\begin{align}
	\sum_{t=1}^{T} Q_t - \lambda \frac{K}{2\tilde{m}} \sum_{t=1}^{T} \nu_{i,t} &\le  36\lambda\frac{K^2}{\tilde{m}^2}B^2\log\left(\delta^{-1}\right)\notag\\ 
	&\le \frac{9}{32} B\log(\delta^{-1})\label{eq:term12}.
	\end{align}
	
	\paragraph{Conclusion:} To conclude, we inject bounds obtain in \eqref{eq:term11} and \eqref{eq:term12} into \eqref{eq:term0}.
	
\end{proof}

We provide a key lemma that will be used in the proof of Theorem~\ref{thm:1} and~\ref{thm:1bis}.
\begin{lemma}\label{lem:keyf1}
	Let $\lambda \in \left(0, \frac{\tilde{m}}{128 K}\bar{\lambda}\right)$, where $\bar{\lambda}$ is defined in \eqref{eq:c_cst}. Consider Algorithm~\ref{algo:1} with inputs $(\lambda , m)$. We have with probability at least $1-9\delta$
	\begin{equation*}
	\sum_{t=1}^{T}\hat{\mu}_t -\frac{7\bar{\lambda}}{32} \sum_{t=1}^{T} \hat{\xi}_t \le \min_{i \in \intr{K}} L_{i,T} + c~\frac{1}{\lambda} \log\left(\frac{\tilde{m}}{B\lambda \delta} \right) 
	\end{equation*}
	where $\tilde{m} = \max\{1,m-1\}$ and $c$ is a numerical constant.
\end{lemma}
\begin{proof}
	\noindent For each $i \in \intr{K}$ and $t\in \intr{T}$, let $\hat{h}_{i,t} := \hat{\ell}_{i,t} - \lambda \hat{v}_{i,t}$ and $h_{i,t} := \mathbb{E}_{t-1}\left[\hat{h}_{i,t}\right]$. Using Lemma~\ref{lem:key0}, we have
	\begin{align}
	\sum_{t=1}^{T} \hat{\mu}_t -\frac{7\bar{\lambda}}{32} \sum_{t=1}^{T} \hat{\xi}_t &\le \min_{i \in \intr{K}}\sum_{t=1}^{T} \hat{h}_{i,t} + \frac{1}{\lambda} \sum_{t=1}^{T} M_t + \frac{\log(K)}{\lambda} + \left(\frac{11 \lambda K}{\tilde{m}} - \frac{7}{32}\bar{\lambda} \right) \sum_{t=1}^{T} \hat{\xi}_t \notag\\
	&\le \min_{i \in \intr{K}}\sum_{t=1}^{T} \hat{h}_{i,t} + \frac{1}{\lambda} \sum_{t=1}^{T} M_t - \frac{\bar{\lambda}}{8} \sum_{t=1}^{T} \hat{\xi}_t + \frac{\log(K)}{\lambda}, \label{bound:TT}
	\end{align}
	where we used the fact that $\lambda \in \left(0, \frac{\tilde{m}}{128K} \bar{\lambda}\right)$.
	
	%		Similarly, if $m=2$ ($\tilde{m}=1$), we plug the bounds obtained in \eqref{eq:term2_1} and \eqref{eq:term2_5} into \eqref{eq:h}
	%		\begin{align}
	%		\sum_{t=1}^{T}\hat{\mu}_t - \frac{3\bar{\lambda}}{32K} \sum_{t=1}^{T} \hat{\xi}_t &\le \min_{i\in \intr{K}} \sum_{t=1}^{T} \hat{h}_{i,t} + \frac{1}{\lambda} \sum_{t=1}^{T}M_t - \frac{\bar{\lambda}}{16K} \sum_{t=1}^{T} \hat{\xi}_t+ \frac{\log(K)}{\lambda} + \left(11\lambda K- \frac{\bar{\lambda}}{32K}\right) \sum_{t=1}^{T} \hat{\xi}_t \notag\\
	%		&\le \min_{i\in \intr{K}} \sum_{t=1}^{T} \hat{h}_{i,t} + \frac{1}{\lambda} \sum_{t=1}^{T}M_t - \frac{\bar{\lambda}}{16K} \sum_{t=1}^{T} \hat{\xi}_t + \frac{\log(K)}{\lambda},\label{bound:TT2}  
	%		\end{align}
	%		where we used the fact that the last term is negative since $\lambda < \frac{\tilde{m}}{352K}\bar{\lambda}$. 
	
	\noindent In order to conclude, we only need bounds on the terms $\min_{i\in \intr{K}} \sum_{t=1}^{T} \hat{h}_{i,t}$ and $\frac{1}{\lambda} \sum_{t=1}^{T} M_t$. Recall that Lemma~\ref{lem:mt} shows that $(M_t)$ is a martingale difference sequence and provides a bound on its conditional variance. Hence, applying Corollary~\ref{cor:concentration_mart2} to this sequence with $c = 3B\bar{\lambda}/40$, with probability at least $1-3\delta$, it holds
	\begin{equation*}
	\frac{1}{\lambda} \sum_{t=1}^{T} M_t - \frac{\tilde{m} \bar{\lambda}}{40 \bar{\lambda}^2 K} \sum_{t=1}^{T} 5\frac{K}{\tilde{m}} \lambda^2 \hat{\xi}_t \le \frac{324K}{\tilde{m}\bar{\lambda}} \left(\log \delta^{-1} + 2 \log_2^{+} \left(\frac{7024}{B^2\bar{\lambda}^2}\right)\right).
	\end{equation*}
	We conclude that 
	\begin{equation}\label{bound:Mt}
	\frac{1}{\lambda} \sum_{t=1}^{T} M_t - \frac{\bar{\lambda}}{8} \sum_{t=1}^{T} \hat{\xi}_t \le 8428 \frac{K}{\tilde{m} \bar{\lambda}} \log\left(\frac{1}{B\bar{\lambda}\delta}\right).
	\end{equation}
	
	%		\noindent Lemma~\ref{lem:Mt} provides the following bound with probability at least $1- 3\delta$
	%		\begin{equation}\label{bound:Mt}
	%		\frac{1}{\lambda} \sum_{t=1}^{T} M_t - \frac{\bar{\lambda}}{8} \sum_{t=1}^{T} \hat{\xi}_t \le 8424~\frac{K}{\tilde{m} \bar{\lambda}}\log\left(\frac{1}{ B \bar{\lambda} \delta}\right).
	%		\end{equation}
	
	\noindent Next, to bound the term $\min_{i\in \intr{K}} \sum_{t=1}^{T} \hat{h}_{i,t}$ we use Lemma~\ref{lem:minL}. We have with probability at least $1-6\delta$
	\begin{align}
	\min_{i\in \intr{K}} \sum_{t=1}^{T} \hat{h}_{i,t} &= \min_{i \in \intr{K}} \hat{L}_{i,T} - \lambda \hat{V}_{i,T}\notag\\
	&\le \min_{i \in \intr{K}} L_{i,T} +\frac{721}{\lambda} \log\left(\frac{\tilde{m}}{B\lambda \delta}\right).\label{bound:minL}
	\end{align}
	
	\noindent Finally, we inject \eqref{bound:Mt} and \eqref{bound:minL} into \eqref{bound:TT} and use $\lambda \in \left(0, \frac{\tilde{m}}{128K} \bar{\lambda}\right)$. We obtain that with probability at least $1-9\delta$
	\begin{equation*}
	\sum_{t=1}^{T}\hat{\mu}_t - \frac{7\bar{\lambda}}{32} \sum_{t=1}^{T} \hat{\xi}_t \le \min_{i \in \intr{K}} L_{i,T} + c~\frac{1}{\lambda} \log\left(\frac{\tilde{m}}{B\lambda \delta} \right),
	\end{equation*}
	where $c$ is a numerical constant.
\end{proof}

The following Lemma is specific to the case $m=p=2$ and $\text{IC}=\text{True}$ in Algorithm~\ref{algo:2}.
\begin{lemma}\label{lem:keyf2}
	Let $\lambda \in \left(0, \frac{\bar{\lambda}}{352 K^2}\right)$, where $\bar{\lambda}$ is defined in \eqref{eq:c_cst}. Consider Algorithm~\ref{algo:2} with input $\lambda $. We have with probability at least $1-9\delta$
	\begin{equation*}
	\sum_{t=1}^{T} \hat{\mu}_t - \frac{3\bar{\lambda}}{32K} \sum_{t=1}^{T} \hat{\xi_t} \le \min_{i \in \intr{K}} L_{i,T} + c~\frac{1}{\lambda} \log\left(\frac{1}{B\lambda \delta} \right),
	\end{equation*}
	where $c$ is a numerical constant.
\end{lemma}
\begin{proof}
	\noindent For each $i \in \intr{K}$ and $t\in \intr{T}$, let $\hat{h}_{i,t} := \hat{\ell}_{i,t} - \lambda \hat{v}_{i,t}$ and $h_{i,t} := \mathbb{E}_{t-1}\left[\hat{h}_{i,t}\right]$. Using Lemma~\ref{lem:key0}, we have
	\begin{align}
	\sum_{t=1}^{T} \hat{\mu}_t -\frac{3\bar{\lambda}}{32K} \sum_{t=1}^{T} \hat{\xi}_t &\le \min_{i \in \intr{K}}\sum_{t=1}^{T} \hat{h}_{i,t} + \frac{1}{\lambda} \sum_{t=1}^{T} M_t + \frac{\log(K)}{\lambda} + \left(11 \lambda K - \frac{3\bar{\lambda}}{32K} \right) \sum_{t=1}^{T} \hat{\xi}_t \notag\\
	&\le \min_{i \in \intr{K}}\sum_{t=1}^{T} \hat{h}_{i,t} + \frac{1}{\lambda} \sum_{t=1}^{T} M_t - \frac{\bar{\lambda}}{16K} \sum_{t=1}^{T} \hat{\xi}_t + \frac{\log(K)}{\lambda}, \label{bound:TT2}
	\end{align}
	where we used the fact that $\lambda \in \left(0, \frac{\bar{\lambda}}{352K^2} \right)$.

	\noindent The remainder of the proof is similar to the proof of Lemma~\ref{lem:keyf1}. 
	
	\noindent Lemma~\ref{lem:mt} provides the following bound with probability at least $1- 3\delta$
	\begin{equation}\label{bound:Mt2}
	\frac{1}{\lambda} \sum_{t=1}^{T} M_t - \frac{\bar{\lambda}}{16K} \sum_{t=1}^{T} \hat{\xi}_t \le \frac{3520}{ \bar{\lambda}}\log\left(\frac{1}{ B \bar{\lambda} \delta}\right).
	\end{equation}
	
	\noindent Moreover, Lemma~\ref{lem:minL} provides the following bound with probability at least $1-6\delta$
	\begin{equation}\label{bound:minL2}
	\min_{i\in \intr{K}} \sum_{t=1}^{T} \hat{h}_{i,t} = \min_{i \in \intr{K}} L_{i,T} +\frac{721}{\lambda} \log\left(\frac{1}{B\lambda \delta}\right). 
	\end{equation}
	
	\noindent Finally, we inject \eqref{bound:Mt2} and \eqref{bound:minL2} into \eqref{bound:TT2}. We obtain that with probability at least $1-9\delta$
	\begin{equation*}
	\sum_{t=1}^{T}\hat{\mu}_t - \frac{3\bar{\lambda}}{32K} \sum_{t=1}^{T} \hat{\xi}_t \le  \min_{i \in \intr{K}} L_{i,T} + c~\frac{1}{\lambda} \log\left(\frac{1}{B\lambda \delta} \right),
	\end{equation*}
	where $c$ is a numerical constant.
\end{proof}

\section{On the sampling strategy in the case $m=p=2, \text{IC}=\text{True}$}
%\textcolor{red}{Question pour GB:} doit-on nous contenir de la preuve pour $m=2$ (ce qui est suffisant pour le papier) ou une preuve pour $m \ge 2$ comme ci-dessous ?
Let $\bm{p}$ denote a distribution over $\intr{K}$. Let $\mathcal{E} = \{A,B\}$ denote a random set of elements in $\intr{K}$, such that $A$ is sampled from $\intr{K}$ following $\bm{p}$ and $B$ is sampled independently and uniformly at random from $\intr{K}$ (possibly $A=B$ and $\cE$ is a singleton).
Therefore, we have for each $u,v \in \intr{K}$, such that $u\neq v$:
\[
\mathbb{P}(\mathcal{E} = \{u,v\}) = \frac{\bm{p}_{u}+\bm{p}_{v}}{K},
\] 
and 
\[
\mathbb{P}(\mathcal{E} = \{u\}) = \frac{\bm{p}_{u}}{K}.
\]

\noindent Finally, let $\bm{p}_{\mathcal{E}}$ denote the restriction of the distribution $\bm{p}$ on $\mathcal{E}$, conditional to $\cE$. Let $X$ denote a random variable following $\bm{p}_{\mathcal{E}}$
\[
\forall i \in \mathcal{E}: \bm{p}_{\mathcal{E}}(X = i) = \bm{p}(X = i| \mathcal{E}) = \frac{\bm{p}_i}{\sum_{j \in \mathcal{E}} \bm{p}_j}.
\]
Let $I$ and $J$ denote two random variables on $\intr{K}$ sampled conditionally to $\cE$, independently following $\bm{p}_{\mathcal{E}}$ (with replacement). 

In this section, we prove two results: the marginal distribution of $I$ on $\intr{K}$ is identical to $\bm{p}$, and a bound on the probabilities of the joint unconditional distribution of $(I,J)$.

\begin{lemma}\label{lem:sample1}
	For each $i \in \intr{K}$,
	\[
	\mathbb{P}\left( I = i\right) = \bm{p}_i.
	\]
\end{lemma}
\begin{proof}
	Fix $i \in \intr{K}$. 
	Let $\mathcal{K}$ denote the set of subsets of $\intr{K}$, constituted of at most two elements.
	
	%For $r \in \{1,2\}$, define $\mathcal{K}(r)$ as the set of $m$-tuples in $\mathcal{K}$,  such that exactly $r$ of its elements are equal to $i$.
	
	For any subset $\bm{a} \in \mathcal{K}$, define
	\[
	\bm{p}_{\bm{a}} := \sum_{i \in \bm{a}} \bm{p}_{i}.
	\]

	\noindent We have
	\begin{align*}
	\mathbb{P}\left(I = i\right) &= \sum_{\bm{a} \in \mathcal{K}} \mathbb{P}\left(I = i, \mathcal{E} = \bm{a}\right)\\
	%&= \mathbb{P}\left(I = i | \mathcal{E} = \{i\}\right)~\mathbb{P}\left( \mathcal{E} = \{i\}\right) + \sum_{\bm{a} \in \mathcal{K}\setminus \{\{i\}\}} \mathbb{P}\left(I = i | \mathcal{E} = \bm{a}\right)~\mathbb{P}\left( \mathcal{E} = \bm{a}\right) \\
	&= \mathbb{P}\left(I = i | \mathcal{E} = \{i\}\right)~\mathbb{P}\left( \mathcal{E} = \{i\}\right) + \sum_{u \in \intr{K}\setminus \{i\}} \mathbb{P}\left(I = i | \mathcal{E} = \{u,i\}\right)~\mathbb{P}\left( \mathcal{E} = \{u,i\}\right)\\
	&= \frac{\bm{p}_{i}}{K} + \sum_{u \in \intr{K}\setminus \{i\}} \frac{\bm{p}_i}{\bm{p}_{u}+\bm{p}_{i}}~\frac{\bm{p}_{u}+\bm{p}_{i}}{K} \\
	&= \frac{\bm{p}_{i}}{K} + \frac{\bm{p}_i}{K} (K-1)\\
	&= \bm{p}_i.
	\end{align*}
	
\end{proof}

\begin{lemma}\label{lem:sample2}
	For each $i,j \in \intr{K}$,
	\[
	\mathbb{P}\left(I = i, J = j\right) \ge \frac{1}{K} \bm{p}_i \bm{p}_j.
	\]
\end{lemma}
\begin{proof}
	Fix $i,j \in \intr{K}$. Let $\mathcal{K}$ denote the set of subsets of $\intr{K}$, constituted of at most two elements. 
	
	Suppose that $i=j$. We have
	\begin{align*}
	\mathbb{P}\left(I=i, J=i \right) &= \sum_{\bm{a} \in \cK} \mathbb{P}\left(I = i, J= i, \mathcal{E} = \bm{a}\right)\\
	&= \sum_{\bm{a} \in \cK} \mathbb{P}\left(I = i, J= i| \mathcal{E} = \bm{a}\right) \mathbb{P}\left(\mathcal{E} = \bm{a}\right)\\
	&= \sum_{\bm{a} \in \cK} \mathbb{P}\left(I = i| \mathcal{E} = \bm{a}\right)^2 \mathbb{P}\left(\mathcal{E} = \bm{a}\right),
	\end{align*}
	where we used the fact that $I$ and $J$ are independent conditionally to $\mathcal{E}$ and that $I$ and $J$ follow the same distribution. We use Jensen's inequality:
	\begin{align*}
	\mathbb{P}\left(I=i, J=i \right) &\ge \left(\sum_{a \in \cK}\mathbb{P}\left(I = i| \mathcal{E} = \bm{a}\right) \mathbb{P}\left(\mathcal{E} = \bm{a}\right)\right)^2\\
	&= \bm{p}_i^2.
	\end{align*}
	Now suppose that $i \neq j$. We have
	\begin{align*}
	\mathbb{P}\left(I=i, J=j \right) &= \mathbb{P}\left(I = i, J= j, \mathcal{E} = \{i,j\}\right)\\
	&= \mathbb{P}\left(I=i| \cE = \{i,j\}\right)\mathbb{P}\left(J=j| \cE = \{i,j\}\right) \mathbb{P}\left(\cE = \{i,j\}\right)\\
	&= \frac{\bm{p}_i}{\bm{p}_i + \bm{p}_j}\frac{\bm{p}_j}{\bm{p}_i + \bm{p}_j} \frac{\bm{p}_i + \bm{p}_j}{K}\\
	&= \frac{\bm{p}_i \bm{p}_j}{K}.
	\end{align*}
\end{proof}

\section{Proof of Theorems~\ref{thm:1} and~\ref{thm:1bis}}

We consider the notation of Algorithms~\ref{algo:1} and~\ref{algo:2}. Let $\hat{\pi}_{ij,t} = \mathbb{P}\left(I_t = i, J_t = j | \mathcal{F}_{t-1} \right)$. Introduce ($\hat{\mu}_t$ and $\hat{\xi}_t$ are the same quantities as in the previous section):
\begin{align*}
\hat{\mu}_t &:= \sum_{i \in \intr{K}} \hat{p}_{i,t} \ell_{i,t},\\
\hat{\nu}_t &:= \frac{1}{2}\sum_{i,j \in \intr{K}} \hat{\pi}_{ij,t}~(\ell_{i,t} - \ell_{j,t})^2\\
\hat{\xi}_t &:= \frac{1}{2}\sum_{i,j \in \intr{K}} \hat{p}_{i,t}\hat{p}_{j,t}~(\ell_{i,t} - \ell_{j,t})^2
%\hat{\epsilon}(\delta, T) &:= \log \delta^{-1} + 2\log\left(1+\log_2\left(\sum_{t=1}^{T} \hat{\nu}_t/B^2\right)_{+}\right).
\end{align*}

\noindent We have, using \eqref{def:Ec} with $c=1/\bar{\lambda}$ (implied by Assumption~\ref{assump}, see Lemma~\ref{lem:assumption1}):
\begin{align*}
\sum_{t=1}^{T} \ell_t\left(\frac{F_{I_t}+F_{J_t}}{2}\right) &\le \sum_{t=1}^{T} \left(\frac{1}{2} \ell_{I_t,t} + \frac{1}{2} \ell_{J_t,t} - \frac{\bar{\lambda}}{2} (\ell_{I_t,t}-\ell_{J_t,t})^2 \right)\\
&= \underbrace{\frac{1}{2} \sum_{t=1}^{T} \mathds{U}_t + \frac{1}{2} \sum_{t=1}^{T} \mathds{U'}_t -\frac{\tilde{m} \bar{\lambda}}{32K} \sum_{t=1}^{T} \hat{\xi}_t - \frac{\bar{\lambda}}{2} \sum_{t=1}^{T} \mathds{W}_t - \frac{\bar{\lambda}}{4} \sum_{t=1}^{T} \hat{\nu}_t}_{\text{Term 1}}\\
&\qquad + \underbrace{\sum_{t=1}^{T} \hat{\mu}_t + \frac{\tilde{m}\bar{\lambda}}{32K} \sum_{t=1}^{T} \hat{\xi}_t - \frac{\bar{\lambda}}{4} \sum_{t=1}^{T} \hat{\nu}_t}_{\text{Term 2}},
\end{align*}
where 
\[
\mathds{U}_t:= \ell_{I_t,t} - \hat{\mu}_t ; \qquad \mathds{U}'_t:= \ell_{J_t,t} - \hat{\mu}_t; \qquad
\mathds{W}_t:= (\ell_{I_t,t} - \ell_{J_t,t})^2 - \hat{\nu}_t .
\]

Section~\ref{sec:BT1} below is common to Theorem~\ref{thm:1} and~\ref{thm:1bis}. In Section~\ref{sec:BT2}, we distinguish between the case where $(p=m=2, \text{IC}=\text{True})$ and $(p=2, m\ge 3)$ or $(p=2, m=2, \text{IC}=\text{False})$.

\subsection{Bounding Term 1}\label{sec:BT1}

Recall that in Algorithm~\ref{algo:1} we have by definition of $I_t$, conditionally to $\cF_{t-1}$: $I_t \sim \hat{p}_t$. Furthermore, in Algorithm~\ref{algo:2}, using Lemma~\ref{lem:sample1}, conditionally to $\cF_{t-1}$, we have: $I_t \sim \hat{p}_t$. Hence,  $(\mathds{U}_t)_{t \in \intr{T}}$ is a martingale difference sequence bounded in absolute value by $B$.
Moreover, we have for all $t \in \intr{T}$
\[
\mathbb{E}\left[\mathds{U}_t^2 | \mathcal{F}_{t-1} \right]  =  \hat{\xi}_t.
\]

\noindent Next we apply the high probability bound provided by Corollary~\ref{cor:concentration_mart2} to the sequence $(\mathds{U}_t)_{t \in \intr{T}}$, with $c= \tilde{m}B\bar{\lambda}/(32K)$. We have with probability at least $1 - 3\delta$
\begin{equation}\label{eq:freed_ut}
\sum_{t=1}^{T} \mathds{U}_t - \frac{\tilde{m}}{32K}\bar{\lambda} \sum_{t=1}^{T} \hat{\xi}_t \le 7700  \frac{K}{\tilde{m}\bar{\lambda}} \log\left(\frac{K}{\tilde{m}B\bar{\lambda}\delta}  \right).
\end{equation}

\noindent Recall that in Algorithm~\ref{algo:1} and~\ref{algo:2}, $I_t$ and $J_t$ have the same marginal distribution. Therefore, with probability at least $1-3\delta$, \eqref{eq:freed_ut} holds with $\mathds{U}_t$ replaced by $\mathds{U}'_t$.

\noindent Similarly, the sequence $((-\bar{\lambda}/2)\mathds{W}_t)_{t \in \intr{T}}$ is a martingale difference bounded in absolute value by $\bar{\lambda}B^2$. For any $t \in \intr{T}$,
\[
\frac{\bar{\lambda}^2}{4}\mathbb{E}\left[\mathds{W}^2_t | \mathcal{F}_{t}\right] \le \frac{\bar{\lambda}^2}{4} \mathbb{E}\left[ \left(\ell_{I_t,t} - \ell_{J_t,t}\right)^4 | \mathcal{F}_{t-1}\right] \le \frac{\bar{\lambda}^2B^2}{4} \hat{\nu}_t.
\] 	

\noindent Next, we apply Corollary~\ref{cor:concentration_mart2} to the sequence $((-\bar{\lambda}/2)\mathds{W}_t)_{t \in \intr{T}}$: We take $c= 1$, we have with probability $1-3\delta$:
\begin{align}
-\frac{\bar{\lambda}}{2} \sum_{t=1}^{T} \mathds{W}_t - \frac{\bar{\lambda}}{4} \sum_{t=1}^{T} \hat{\nu}_t &\le 72 \bar{\lambda}B^2 \log(\delta^{-1}) \notag \\
&\le 72 B \log(\delta^{-1}).\label{eq:freed_wt}
\end{align}

\noindent Using \eqref{eq:freed_ut} and \eqref{eq:freed_wt}, we conclude that with probability $1- 9 \delta$

\begin{equation}\label{eq:term_1}
\text{Term 1} \le 7772 \frac{K}{\tilde{m}\bar{\lambda}} \log\left(\frac{K}{\tilde{m}B\bar{\lambda}\delta}  \right).
\end{equation}

\subsection{Bounding Term 2}\label{sec:BT2}
We divide this part of the proof into two section (depending on the expression of the joint distribution $\hat{\pi}_t$).

\subsubsection{Case  ($p=2$ and $m \ge 3$) or ($p=2$, $m = 2$ and $\text{IC}=\text{False}$)}

Recall that conditionally to $\cF_{t-1}$, the played experts $I_t$ and $J_t$ are sampled independently according to $\hat{p}_t$ from $\intr{K}$. Therefore for any $i,j \in \intr{K}$, $\hat{\pi}_{ij,t} = \hat{p}_{i,t}\hat{p}_{j,t}$ and $\hat{\nu}_t = \hat{\xi}_t$.

\noindent Hence, Term 2 satisfies the following bound
\begin{equation*}
\text{Term 2} \le \sum_{t=1}^{T} \hat{\mu}_t - \frac{7\bar{\lambda}}{32} \sum_{t=1}^{T} \hat{\xi_t}.
\end{equation*}

\noindent Using the first claim of Lemma~\ref{lem:keyf1}, we have if $\lambda \in \left(0, \frac{\tilde{m}}{128K} \bar{\lambda}\right)$
\begin{equation}\label{eq:term_2}
\text{Term 2} \le \min_{i \in \intr{K}} L_{i,T} + c~\frac{1}{\lambda} \log\left(\frac{\tilde{m}}{B\lambda \delta} \right),
\end{equation}
where $c$ is a numerical constant. The conclusion of the theorem follows by combining the upper bounds obtained in \eqref{eq:term_1} and \eqref{eq:term_2}.

\subsubsection{Case $m=p=2$ and $\text{IC}=\text{True}$:}

Using Lemma~\ref{lem:sample1} we have $I_t \sim \hat{p}_t$. Furthermore, using Lemma~\ref{lem:sample2} we have that for any $i,j \in \intr{K}$, any $t \in \intr{T}$:
\[
\hat{\pi}_{ij,t} \ge \frac{1}{K} \hat{p}_{i,t} \hat{p}_{j,t}.
\]
Therefore $\hat{\nu}_t \ge \frac{1}{K} \hat{\xi}_t$, and we have the following bound on Term 2:
\begin{equation*}
\text{Term 2} \le \sum_{t=1}^{T} \hat{\mu}_t - \frac{3\bar{\lambda}}{32K} \sum_{t=1}^{T} \hat{\xi_t}.
\end{equation*}
\noindent Using the second claim of Lemma~\ref{lem:keyf2}, we have if $\lambda \in \left(0, \frac{\bar{\lambda}}{352K^2}\right)$
\begin{equation}\label{eq:term_2b}
\sum_{t=1}^{T} \hat{\mu}_t - \frac{7}{32B} \sum_{t=1}^{T} \hat{\xi_t} \le \min_{i \in \intr{K}} L_{i,T} + c~\frac{1}{\lambda} \log\left(\frac{1}{\lambda B \delta} \right).
\end{equation}
The conclusion of the theorem follows by combining the upper bounds obtained in \eqref{eq:term_1} and \eqref{eq:term_2b}.

\section{Proofs of lower bounds, Theorem~\ref{thm:2} and Theorem~\ref{thm:4}}
The proofs of Theorem~\ref{thm:2} and Theorem~\ref{thm:4} are presented in four steps. The only difference between the proofs is in the last step. Thus the first three steps are common to both proofs.

We adapt the main steps of \cite{auer1995gambling} to our setting.  The gist of the proof is the following. We construct a distribution with very correlated experts. In this situation, going from a weighted average of experts to a single expert with the largest weight does not change the prediction risk much. Then, we use some classical arguments in deriving lower bounds for the expected regret using information theory results.

Let $T>0$ be fixed, we consider that the loss function is the squared loss and we focus on the particular setting where the target variables $(Y_t)$ are identically $0$.
\paragraph{First step: Specifying the distributions.}
We start by considering a deterministic forecaster. We denote by $\mathbb{P}_i$ the joint distribution of expert predictions, where all experts are identical and distributed as one and the same Bernoulli variable with parameter $1/2$, except the optimal expert $i$ who has distribution $\mathcal{B}\left( \frac{1}{2} - \epsilon \right)$ but
is still strongly correlated to the others.

More precisely, let $(U_t)_{t \in \intr{T}}$ be a sequence of independent random variables distributed according the uniform distribution on $[0,1]$. We consider that in each round the expert predictions have the following joint distribution $\mathbb{P}_i$:
\begin{itemize}
	\item For $j \neq i$: $F_{j,t} = \mathds{1}\left(U_t \le \frac{1}{2}\right)$.
	\item $F_{i,t} = \mathds{1}\left( U_t \le \frac{1}{2} - \epsilon\right)$.
\end{itemize}
Recall that in this setting we have for any $k,j \in \intr{K}\setminus \{i\}$
\begin{align*}
\mathbb{E}_i[F_{j,t}F_{k,t}] &= \frac{1}{2}\\
\mathbb{E}_i[F_{i,t}F_{j,t}] &= \frac{1}{2} - \epsilon.
\end{align*}
Finally, we denote by $\mathbb{P}_0$ the joint distribution where all experts are equal to the
same Bernoulli$(1/2)$ variables, i.e., experts predictions are defined by $F_{i,t} = \mathds{1}(U_t \le 1/2)$,
$i\in \intr{K}$.

\paragraph{Second step: Strategy Reduction.}
Suppose that the player follows a deterministic strategy $\mathcal{A}$. In each round $t$, given $\mathcal{F}_{t-1}$, this strategy selects a subsets $S_t$ of $\intr{K}$ of size $m$ and a sequence of non-negative weights $(\alpha_{i,t})_{i \in S_t}$, such that $\sum_{i} \alpha_{i,t} = 1$, and plays the convex combination $\sum_{i \in S_t} \alpha_{i,t} F_{i,t}$. % then observes all the experts predictions $(F_{i,t})_{i \in S_t}$.

For such a strategy $\mathcal{A}$, we associate a strategy $\hat{\mathcal{A}}$, such that in each round, we run the strategy $\mathcal{A}$ except that we play only the expert with the largest weight $\hat{i}_t \in \argmax_{i \in S_t} \alpha_{i,t}$.

Let us analyse the difference of the cumulative loss between the strategies $\mathcal{A}$ and $\hat{\mathcal{A}}$. Let $l_t(\mathcal{A})$ denote the loss of the strategy $\mathcal{A}$ at round $t$. We have

\begin{equation*}
\mathbb{E}_i\left[l_t\left(\mathcal{A}\right) - l_t(\hat{\mathcal{A}})\right] = \ee{i}{\paren[3]{\sum_{j \in S_t} \alpha_{j,t} F_{j,t}}^2} - \ee{i}{\paren[3]{ \sum_{j \in S_t} \mathds{1}\paren{\hat{i}_t = j}F_{j,t} }^2 }.
\end{equation*}

If $i \notin S_t$ then we have $\mathbb{E}_i[l_t\left(\mathcal{A}\right) - l_t(\hat{\mathcal{A}})] = 0$. 

If $i \in S_t$ and $\hat{i}_t = i$, we have (let $j \in \intr{K}$ such that $j \neq i$)
\begin{align*}
\mathbb{E}_i\left[l_t\left(\mathcal{A}\right) - l_t(\hat{\mathcal{A}})\right]
&= \mathbb{E}_i\left[ \left((1-\alpha_{i,t})F_{j,t} + \alpha_{i,t} F_{i,t} \right)^2\right] - \mathbb{E}_i\left[F_{i,t}\right]\\
&= %\mathbb{E}_{i}\left[
\left(1 - \alpha_{i,t}\right)^2 \frac{1}{2} + \alpha_{i,t}^2 \left(\frac{1}{2}- \epsilon\right) + 2\alpha_{i,t}(1-\alpha_{i,t}) \left(\frac{1}{2} - \epsilon\right)%\right]
- \frac{1}{2}+\epsilon\\
&=%\mathbb{E}_i \left[
\epsilon \left(1 - \alpha_{i,t}\right)^2 %\right]
\\
&\ge 0.
\end{align*}	
If $i \in S_t$ and $\hat{i}_t \neq i$, we have (let $j \in \intr{K}$ such that $j \neq i$) 	
\begin{align*}
\mathbb{E}_i\left[l_t\left(\mathcal{A}\right) - l_t(\hat{\mathcal{A}})\right]
&= \mathbb{E}_i \left[ \left((1-\alpha_{i,t})F_{j,t} + \alpha_{i,t} F_{i,t} \right)^2\right] - \mathbb{E}_i \left[F_{j,t}\right]\\
&= %\mathbb{E}_i \left[
\left(1 - \alpha_{i,t}\right)^2 \frac{1}{2} + \alpha_{i,t}^2 \left(\frac{1}{2}- \epsilon\right) + 2\alpha_{i,t}(1-\alpha_{i,t}) \left(\frac{1}{2} - \epsilon\right) %\right]
- \frac{1}{2}\\
&= %\mathbb{E}_i \left[
\epsilon \alpha_{i,t}^2 - 2\epsilon \alpha_{i,t}%\right]
\\
&\ge -\frac{3}{4} \epsilon,
\end{align*}
where we used the fact that $\alpha_{i,t} \in [0,1/2]$, since $\hat{i}_t \neq i$.

To summarize, in the worst case, the excess loss between $\mathcal{A}$ and $\hat{\mathcal{A}}$ is $-\frac{3}{4} \epsilon$. Hence, we have the following lower bound on the expected regret between the two strategies:
\begin{equation}\label{eq:reduction}
\mathcal{R}_T (\mathcal{A}) - \mathcal{R}_T (\hat{\mathcal{A}}) \ge -\frac{3}{4}T\epsilon.
\end{equation}

\paragraph{Third step: Information theoretic tools.}

Let us introduce the following notation: assume the player follows a deterministic strategy $\cA$, and let
$Z_t=(C_t,\l_t(F_{i,t})_{i \in C_t})$ denote the information disclosed to the player at time $t$.
Denote $\bZ^t=(Z_1,\ldots,Z_t)$ the entire information available to the player since the start.
The quantities $Z_t,\bZ^t$ are considered as random variables, whose distribution is determined by the
underlying experts distribution, and the player strategy $\cA$. 

\begin{lemma}\label{lem:info}
	Let $F(\bZ^T)$ be any fixed function of the player observations, taking values in $[0,B]$.
	Then for any $i \in \intr{K}$ and any player strategy $\cA$,
	\begin{equation*}
	\ee{i}{F\paren{\bZ^T}} \le \ee{0}{F\paren{\bZ^T}} +\frac{B}{2} \sqrt{\ee{0}{N_i} \log(1-2\epsilon)^{-1}},
	\end{equation*}
	where $N_i = \sum_{i=1}^T \ind{i \in C_t}$.
	
	In the case where $\abs{C_t}=1$ for all $t$, the following sharper bound holds:
	\begin{equation*}
	\ee{i}{F\paren{\bZ^T}} \le \ee{0}{F\paren{\bZ^T}} +\frac{B}{2} \sqrt{\ee{0}{N_i} \log(1-4\epsilon^2)^{-1}},
	\end{equation*}
	
\end{lemma}

\begin{proof}
	Fix $i\in \intr{K}$.  Denote $\mbq_i$ the distribution
	of $\bZ^T$ induced by expert distribution $\mbp_i$ and a fixed player strategy $\cA$ (omitted from the notation
	for simplicity). For any function $G$ bounded by $R$, it is well-known that it holds $\abs{\ee{X\sim \mbp}{G(X)} -
		\ee{X\sim\mbq}{G(X)}} \leq 2R \norm{\mbp-\mbq}_{TV}$, where $\norm{\cdot}_{TV}$ denotes the total variation distance.
	Hence, by shifting $F$ by $-B/2$, we get
	\[
	\ee{i}{F(\bZ^T)} - \ee{0}{F(\bZ^T)} \leq B \norm{\mbq_i-\mbq_0}_{TV} \leq B \sqrt{\frac{1}{2} 
		\text{KL}\paren{\mbq_0 \| \mbq_i}},
	\]
	by Pinsker's inequality, where $\text{KL(.)}$ denotes the Kullback-Leibler divergence.

	Next, we will compute the quantity $\text{KL}\paren{\mbq_0 \| \mbq_i}$.
	The chain rule for relative entropy (Theorem 2.5.3 in \citealp{cover1999elements}) gives:
	\begin{equation} \label{eq:chainkl}
	\text{KL}\paren{\mbq_0 \| \mbq_i}
	= \sum_{t=1}^{T} \text{KL}\paren{\mbq_0\set{Z_t | \bZ^{t-1}} \| \mbq_i\set{Z_t | \bZ^{t-1}}},
	\end{equation}
	where 
	\begin{align*}
	\text{KL}\paren{\mbq_0\set{Z_t | \bZ^{t-1}\} \| \mbq_i\{Z_t | \bZ^{t-1}}} & :=
	\sum_{\bz^{t}}
	\mbq_0\set{\bz^{t-1}}\mbq_0\set{z_t|\bz^{t-1}} \log\paren{\frac{\mbq_0\set{z_t|\bz^{t-1}}}{\mbq_i\set{z_t|\bz^{t-1}}}}\\
	& = \sum_{\substack{\bz^{t}\\\text{s.t. } i \in C_t\\}}
	\mbq_0\set{\bz^{t-1},C_t} %\mbq_0\set{C_t|\bz^{t-1}}
	\mbq_0\set{z_t|C_t} \log\paren{\frac{\mbq_0\set{z_t|C_t}}{\mbq_i\set{z_t|C_t}}}.
	% + \mbq_0\set{i \not\in C_t} \sum_{\bz^{t-1}}
	% \mbq_0\set{\bz^{t-1} | i \not\in C_t} \log\paren{\frac{\mbq_0\set{z_t|\bz^{t-1}}}{\mbq_i\set{z_t|\bz^{t-1}}}}.
	\end{align*}
	The last line holds because $\mbq_\bullet\set{z_t|\bz^{t-1}}=\mbq_\bullet\set{z_t|\bz^{t-1},C_t}
	\mbq_\bullet\set{C_t|\bz^{t-1}}$, and it holds $\mbq_0\set{C_t|\bz^{t-1}} = \mbq_i\set{C_t|\bz^{t-1}}$
	since the strategy's play only depends on past observations; also $\mbq_\bullet\set{z_t|\bz^{t-1},C_t} = \mbq_\bullet \set{z_t|C_t}$
	since the observed experts' losses at round $t$ are independent of the past given the choice of $C_t$.
	Furthermore, if $i\not\in C_t$, one has $\mbq_0\set{z_t|C_t}=\mbq_i\set{z_t|C_t}$.

	On the other hand, if $z_t$ is such that $i \in C_t$, then:
	\begin{itemize}
		\item under $\mbq_0$ since all experts are identical and equal
		to the same Ber$(1/2)$ variable (and $Y_t$ is identically 0), $\mbq_0(z_t|C_t)$ only charges the two points with all
		observed losses equal to 0 (denote this $u_0$) or all equal to 1 (denote this $u_1$), each with probability $1/2$;
		\item under $\mbq_i$, it holds $\mbq_i(u_1|C_t) =\frac{1}{2}-\epsilon$ and $\mbq_i(u_0|C_t) \geq \frac{1}{2}$.
		In fact, if $\abs{C_t}\geq 2$, then $\mbq_i(u_0|C_t)=\frac{1}{2}$ (since with probability $\epsilon$ under
		$\mbq_i$, we observe a state that is neither $u_0$ nor $u_1$, namely when all observed experts err but $F_i$),
		and if $\abs{C_t}=1$, then $\mbq_i(u_0|C_t)=\frac{1}{2} + \eps$ (since $F_i$ alone is observed then).
	\end{itemize} 
	
	Therefore, in general
	\begin{align*}
	\text{KL}\paren{\mbq_0\set{Z_t | \bZ^{t-1}\} \| \mbq_i\{Z_t | \bZ^{t-1}}}
	%	\text{KL}\left(\mathbb{P}_0\{l_t | \mathbf{l}^{t-1}, i \in C_t\} \| \mathbb{P}_i\{l_t | \mathbf{l}^{t-1}, i \in C_t\}\right)
	&\leq \mbp_0(i \in C_t ) \paren{\frac{1}{2}\log\paren{\frac{1/2}{1/2 - \epsilon}} +\frac{1}{2}\log\paren{\frac{1/2}{1/2}}}\\
	&\leq \frac{1}{2}\mbp_0(i \in C_t ) {\log\paren{1 - 2 \epsilon}^{-1}} . 
	\end{align*}
	In the case where $\abs{C_t}=1$ for all $t$, we get the sharper bound
	\begin{align*}
	\text{KL}\paren{\mbq_0\set{Z_t | \bZ^{t-1}\} \| \mbq_i\{Z_t | \bZ^{t-1}}}
	&= \mbp_0(i \in C_t ) \paren{\frac{1}{2}\log\paren{\frac{1/2}{1/2 - \epsilon}} +\frac{1}{2}\log\paren{\frac{1/2}{1/2+\epsilon}}}\\
	&= \frac{1}{2}\mbp_0(i \in C_t ) {\log\paren{1 - 4 \epsilon^2}^{-1}} . 
	\end{align*}
	Plugging this into~\eqref{eq:chainkl}, we obtain 
	
	$	\text{KL} \left(\mbq_0 \| \mbq_i\right) \leq -\frac{1}{2} \mathbb{E}_0\left[N_i\right] \log\left(1 - 2\epsilon\right),$ resp. $	\text{KL} \left(\mbq_0 \| \mbq_i\right) \leq -\frac{1}{2} \mathbb{E}_0\left[N_i\right] \log\left(1 - 4\epsilon^2\right),$ if $\abs{C_t}=1$ for all $t$, leading to the claims.
	
\end{proof}

\paragraph{Fourth step for Theorem~\ref{thm:2}: lower bounding the regret of $\hat{\mathcal{A}}$ in the case $\abs{C_t}\ge 2$.}
Recall $\hat{i}_t$ denotes the single expert played by the ``reduced'' strategy $\hat{\cA}$.
At round $t$, the expected loss for the player playing $\hat{\cA}$ is given by
\begin{align*}
\mathbb{E}_i\left[l_{t,\hat{i}_t}\right] &= \left(\frac{1}{2} - \epsilon\right) \mathbb{P}_i\left(\hat{i}_t = i\right) + \frac{1}{2} \mathbb{P}_i\left(\hat{i}_t \neq i\right)
= \frac{1}{2} - \epsilon~\mathbb{P}_i\left(\hat{i}_t = i\right).
\end{align*}
For each $j \in \intr{K}$ let $M_j := \sum_{t=1}^{T} \mathds{1}\set{\hat{i}_t = j}$. Hence
\begin{equation*}
\sum_{t=1}^{T} \mathbb{E}_i \left[l_{t,\hat{i}_t}\right] = \frac{T}{2} - \epsilon~\mathbb{E}_i \left[M_i\right],
\end{equation*}
and the regret with respect to the optimal arm $i$ under $\mathbb{P}_i$ is 
\begin{equation}
\label{eq:regri}
\mathbb{E}_i \left[\mathcal{R}_T(\hat{\mathcal{A}})\right] = \epsilon \left(T - \mathbb{E}_i \left[M_i\right] \right).
\end{equation}

We can apply Lemma~\ref{lem:info} to $F(\bZ^t)=M_i$:
since we assume the player follows a deterministic
strategy, $M_i$ is a function of the information $\bZ^t$ available to the player, bounded by $T$. Thus it holds:
\begin{equation}
\mathbb{E}_i\left[M_i\right] \le \mathbb{E}_0 \left[M_i\right] + \frac{T}{2} \sqrt{\mathbb{E}_0\left[N_i\right] \log\left(1 - 2\epsilon\right)^{-1} }.
\end{equation}
Observe that $\sum_{i=1}^K M_i = T$ and $\sum_{i=1}^K N_i = mT$. Hence
\begin{align*}
\sum_{i=1}^K	\mathbb{E}_i\left[M_i\right] &\le \sum_{i=1}^{K} \mathbb{E}_0\left[M_i\right] + \frac{T}{2} \sum_{i=1}^{K} \sqrt{\mathbb{E}_0\left[N_i\right] \log\left(1 - 2\epsilon\right)^{-1}} \\
&\le \mathbb{E}_0\left[\sum_{i=1}^{K} M_i \right] + \frac{TK}{2}\sqrt{ \frac{1}{K} \sum_{i=1}^{K} \mathbb{E}_0\left[  N_i\right] \log\left(1 - 2\epsilon\right)^{-1}}\\
&  = T + {T^{\frac{3}{2}} \sqrt{mK \epsilon  }},
\end{align*}
where we used the fact that for $\epsilon \in \left(0, 1/4\right):  -\log\left(1 - 2\epsilon\right) \le 4\epsilon$. Let $\mbp_*=\frac{1}{K} \sum_{i=1}^K \mbp_i$ the adversary choosing uniformly at random
among the expert distributions $\mbp_i$ at the start of the game (i.e. choosing at random the optimal expert).
From the above and~\eqref{eq:regri} we deduce
\begin{align*}
\ee{*}{\mathcal{R}_T(\hat{\mathcal{A}})} \geq \frac{1}{K} \sum_{i=1}^K
\ee{i}{\mathcal{R}_T(\hat{\mathcal{A}})} 
\geq \epsilon \paren{ T\paren{1-\frac{1}{K}} - T^{\frac{3}{2}} \sqrt{\frac{m \epsilon}{K}}}
\end{align*}

Using inequality \eqref{eq:reduction}, we obtain
\begin{equation*}
\ee{*}{\mathcal{R}_T({\mathcal{A}})}
\geq \epsilon \paren{ T\paren{\frac{1}{4}-\frac{1}{K}} - T^{\frac{3}{2}} \sqrt{\frac{m \epsilon}{K}}}
\geq \epsilon T\paren{ \frac{1}{20} - \sqrt{\frac{Tm\epsilon}{K}}},
%  \ge \epsilon \left(\frac{T}{4} - \frac{T}{K} - T \sqrt{\log(2) \frac{m}{K} \epsilon T}  \right).
\end{equation*}
if $K \geq 5$. Choosing $\epsilon = \frac{1}{900}\frac{K}{mT}$, we get
\begin{equation*}
\ee{*}{\mathcal{R}_T(\cA)} \ge 10^{-5}~\frac{K}{m}.
\end{equation*}

Recall that this lower bound was derived for deterministic players. Generalizing this bound to random players follows simply by applying Fubini's theorem. Also since the bound is in expectation over
expert predictions drawn according to $\mbp_*$, for any strategy $\cA$ there exists at least one deterministic sequence of
expert forecasts with regret larger than its expectation. 

\paragraph{Fourth step for Theorem~\ref{thm:4}: lower bounding the regret of $\hat{\mathcal{A}}$ in the case $\abs{C_t} = 1$.}
The only difference between the proof in this case and the proof in the previous case is the bound given by Lemma~\ref{lem:info}. The regret with respect to the optimal arm $i$ under $\mathbb{P}_i$ is
\begin{equation}\label{eq:reg_m1}
\mathbb{E}_i\left[ \mathcal{R}_T(\hat{\mathcal{A}})\right] = \epsilon (T - \mathbb{E}_i[M_i]).
\end{equation}

We can apply Lemma~\ref{lem:info} to $F(\bZ^t)=M_i$:
since we assume the player follows a deterministic
strategy, $M_i$ is a function of the information $\bZ^t$ available to the player, bounded by $T$. Thus it holds:
\begin{equation*}
\mathbb{E}_i\left[M_i\right] \le \mathbb{E}_0 \left[M_i\right] + \frac{T}{2} \sqrt{\mathbb{E}_0\left[N_i\right] \log\left(1 - 4\epsilon^2\right)^{-1} }.
\end{equation*}
Observe that $\sum_{i=1}^K M_i = T$ and $\sum_{i=1}^K N_i = T$. Hence
\begin{align*}
\sum_{i=1}^K	\mathbb{E}_i\left[M_i\right] &\le \sum_{i=1}^{K} \mathbb{E}_0\left[M_i\right] + \frac{T}{2} \sum_{i=1}^{K} \sqrt{\mathbb{E}_0\left[N_i\right] \log\left(1 - 4\epsilon^2\right)^{-1}} \\
&\le \mathbb{E}_0\left[\sum_{i=1}^{K} M_i \right] + \frac{TK}{2}\sqrt{ \frac{1}{K} \sum_{i=1}^{K} \mathbb{E}_0\left[  N_i\right] \log\left(1 - 2\epsilon^2\right)^{-1}}\\
&  = T + {T^{\frac{3}{2}} \sqrt{2K \epsilon^2  }},
\end{align*}
where we used the fact that for $\epsilon \in \left(0, 1/4\right):  -\log\left(1 - 4\epsilon^2\right) \le 8\epsilon^2$. Let $\mbp_*=\frac{1}{K} \sum_{i=1}^K \mbp_i$ the adversary choosing uniformly at random
among the expert distributions $\mbp_i$ at the start of the game (i.e. choosing at random the optimal expert).
From the above and~\eqref{eq:reg_m1} we deduce
\begin{align*}
\ee{*}{\mathcal{R}_T(\hat{\mathcal{A}})} \geq \frac{1}{K} \sum_{i=1}^K
\ee{i}{\mathcal{R}_T(\hat{\mathcal{A}})} 
\geq \epsilon \paren{ T\paren{1-\frac{1}{K}} - T^{\frac{3}{2}} \sqrt{2\frac{ \epsilon^2}{K}}}
\end{align*}

Using inequality \eqref{eq:reduction}, we obtain
\begin{equation*}
\ee{*}{\mathcal{R}_T({\mathcal{A}})}
\geq \epsilon \paren{ T\paren{\frac{1}{4}-\frac{1}{K}} - T^{\frac{3}{2}} \sqrt{2\frac{ \epsilon^2}{K}}}
\geq \epsilon T\paren{ \frac{1}{20} - \sqrt{2\frac{T\epsilon^2}{K}}},
\end{equation*}
if $K \geq 5$. Choosing $\epsilon = \frac{1}{30}\sqrt{\frac{K}{T}}$, we get
\begin{equation*}
\ee{*}{\mathcal{R}_T(\cA)} \ge 10^{-5}~\sqrt{KT}.
\end{equation*}

The generalization for the random players follows directly using the same argument as in the fourth step of the proof of Theorem~\ref{thm:2}.

\section{Proof of Theorem~\ref{thm:5}}

Let $\ell$ be the squared loss:  $l(x,y)=(x-y)^2$ on $\cX=\cY=[0,1]$.
Consider the game protocol presented in Algorithm~\ref{algo:gp} with $p=1$ and $m\in \intr{K}$. Suppose that the target variable $y$ is identically equal to $0$ ($y_t = 0$ for all $t \in \intr{T}$).
Suppose that at each round $t \in \intr{T}$, for each expert $i \in \intr{K}$, the prediction $F_{i,t}$ follows a Bernoulli distribution of a parameter denoted $\ell_{i,t}$. We have
\begin{equation*}
\mathbb{E}\left[\mathcal{R}_T\right] = \sum_{t=1}^{T} \mathbb{E}\left[F_{I_t,t}\right] - \min_{i \in \intr{K}} \sum_{t=1}^{T} \mathbb{E}\left[F_{i,t}\right].
\end{equation*}
The game protocol presented in Algorithm~\ref{algo:gp} reduces to the $K$-armed bandit game with $m$ feedbacks in each round, analysed in \cite{seldin2014prediction}.

Theorem below presented in \cite{seldin2014prediction} (the full version including appendices) as Theorem 2, provides a lower bound for the regret.
\begin{theorem}[\cite{seldin2014prediction}]\label{thm:seldin}
	For the $K$-armed bandit game with $mT$ observed rewards and $T \ge \frac{3}{16} \frac{K}{m}$,
	\[
	\inf \sup \mathbb{E}\left[\mathcal{R}_T\right] \ge 0.03 \sqrt{\frac{K}{m}T},
	\]
	where the infinimum is over all playing strategies and the supremum is over all individual sequences.
\end{theorem}

The result stated in Theorem~\ref{thm:5} is a direct consequence of the Theorem~\ref{thm:seldin} and the setting described above.

\section{Some implementation details and algorithmic complexity}
\label{app:algocomplex}

We discuss here some details of the implementation of Algorithms~\ref{algo:0},~\ref{algo:1},~\ref{algo:2},
more specifically concerning the cost of keeping track of the distribution $\hat{p}_t$ and of sampling from it at
each round. We concentrate on Algorithm~\ref{algo:1} for simplicity, but the arguments below apply to all
algorithms.

We start with a fundamental observation. While the definitions~\eqref{eq:deflhat},~\eqref{eq:defvhat}
for $\hat{\ell}_{i,t}$ and $\hat{v}_{i,t}$ were
written in order to emphasize the unbiased character of the loss estimates,
the algorithm is unchanged if we use instead the shifted ``pseudo-loss'' estimates
\begin{align}\label{eq:defltilde}
\tilde{\ell}_{i,t}
&  := \hat{\ell}_{i,t}- \ell_{I_t,t} = 
\frac{K}{\tilde{m}}\mathds{1}\left(i \in \mathcal{U}_t\right)(\ell_{i,t}-\ell_{I_t,t}),
\end{align}
and further observe that it holds  $\hat{v}_{i,t} = \tilde{\ell}_{i,t} ^2.$ Using the above pseudo-losses
in place of the estimated losses does not change the sampling distribution $\hat{p}_t$, since
all estimated losses are shifted by the {\em same} quantity $\ell_{I_t,t}$, which gets cancelled
through the normalization in the definition~\eqref{eq:defphat} of the EW distribution $\wh{p}_t$.

Observe that the pseudo-loss estimates $\tilde{\ell}_{i,t}$ (as well as the corresponding
variance estimates $\hat{v}_{i,t}$) are equal to zero for all $i \not\in \cU_t$. Therefore,
to keep track of the cumulative pseudo-loss estimates $\tilde{L}_{i,t} = \sum_{k\leq t} \tilde{\ell}_{i,k}$, only $\abs{\cU_t}=\max\{m-2,1\}$ of them
have to be updated at each round.

In order to keep track and sample efficiently from $\hat{p}_t$, we propose the following construction.
Let $T$ be a balanced binary tree of depth $\lceil \log_2 (K) \rceil$,
with $K$ leaves, such that each leaf $i \in \partial T$ is identified to an expert index.
Furthermore, assume that each internal node $u$ of $T$ stores the partial sum
$S_{u,t} = \sum_{v \in \partial T_{u}}\exp\paren{-\lambda \tilde{L}_{v,t} + \lambda^2 \hat{V}_{v,t}}$, where $T_{u}$
is the subtree of $T$ rooted at node $u$. Then, by the above considerations, it holds that
$S_{u,t} = D_t \sum_{v \in \partial T_u} \hat{p}_{u,t} = D_t \hat{p}_t(\partial T_u)$, where $D_t$ is a factor depending only on $t$
but not on the node $u$. Note also that $D_t = S_{\emptyset}$, where $\emptyset$ denotes the root note of $T$.
It is then possible to sample efficiently $I_t \sim \hat{p}_t$ in a standard manner, as follows:
\begin{enumerate}
	\item Generate $U \sim \mathrm{Unif}[0,1]$, and put $Z=S_{\emptyset} U$. Let $v=\emptyset$.
	\item If $v$ is a leaf of $T$, stop and output $v$.
	\item Let $v_{\mathrm{left}}, v_{\mathrm{right}}$ denote the two descendent nodes of $v$.
	\item If $Z<S_{v_{\mathrm{left}}}$, then let $v \leftarrow v_{\mathrm{left}}$ and go to step 2.
	\item Otherwise, i.e. $Z\geq S_{v_{\mathrm{left}}}$, let $v \leftarrow v_{\mathrm{right}}$,
	$Z \leftarrow Z - S_{v_{\mathrm{left}}}$, and go to step 2.
\end{enumerate}
It easy to check that the above sampling returns a random sample from the probability $\hat{p}_t$.
(Namely, each time that step 2 is reached, conditionally to past steps $Z$ is uniformly distributed in the interval $[0,S_v]$,
and therefore the left or right descendent of $u$ is picked with probability $\hat{p}_t(\partial T_{v_{\mathrm{left}}} | \partial T_v)$ resp. $\hat{p}_t(\partial T_{v_{\mathrm{right}}} | \partial T_v)$; the chain rule yields the claim.)
Obviously, the computing complexity of the above is $\cO(\log K)$ (the depth of the tree).

Furthermore, to update the quantities stored at the nodes of $T$ at each round, since only the estimated
cumulative pseudo-losses of experts $i \in \cU_t$ have their value modified, it is sufficient to do
the following for each $i \in \cU_t$:
\begin{enumerate}
	\item Let $v$ be the leaf representing $i$. Update $S_{v} \leftarrow S_{v} \exp\paren{-\lambda \tilde{\ell}_{i,t} + \lambda^2 \hat{v}_{i,t}}$.
	\item Go up the tree to the root and sequentially update all ancestors $w$ of $v$ according to $S_{w}=S_{w_{\mathrm{left}}} +  S_{w_{\mathrm{right}}}$.
\end{enumerate}
Again, the computing complexity of this update operation is $\cO(\log K)$.

All in all, the computational cost of the initialization of the tree is $\cO(K)$, but then at each round the
computational cost of the sampling and update operations is $\cO(m\log(K))$.

	%\newpage
	
	%\bibliographystyle{abbrv}
	%\bibliography{bib_data_base}
	
	%\vskip 0.2in
	%\bibliographystyle{plain}
	%\bibliography{oomp}
	%\input{oomp.bbl}
\end{document}